\documentclass[12pt]{article}
\usepackage{CJK}
\begin{document}
\baselineskip 5mm
\def\bibite#1#2{\centerline{\hbox to0.7cm{#1\hss}
\parbox[t]{13cm}{#2}}\vspace{\baselineskip}}

\newcommand{\bm}[1]{\mbox{\boldmath $#1$}}
\renewcommand{\theequation}{\thesection.\arabic{equation}}


 \title{{ $\Gamma-\Sigma-C^{0}-$determinacy of $\Gamma-\Sigma-C^{0}-$equivariant bifurcation problems with respect to $\Gamma-\Sigma-C^{0}-$BD and contact equivalence from the weighted point view}
 \thanks{ This work was supported by the National Nature Science Foundation of China under grant, No.10671009, No.60534080}}
 \author{  Suhui Liu\\
   School of Science , Wuhan Institute of Technology\\ Wuhan, Hubei, P.R. of China.\\
   E-mail:17120801@wit.edu.cn,\\
   \\
   Liu~Hengxing\\
   School of Mathematics and Statistics,  Wuhan University,\\Wuhan, 430072. Hubei, P.R. of China\\
   E-mail: jwluan@whu.edu.cn.\\ }
 \date{}
 \maketitle
 \begin{center}
 \begin{minipage}{135mm}
 \begin{center} {\bf abstract}\end{center}
          {\ \ \  In this paper, $C^{0}$ finite determination of $\Gamma-$equivariant bifurcation problems in the relative case from the weighted point view is being discussed . Some criteria on the $C^{0}$ finite determination of $\Gamma-$equivariant bifurcation problems in the relative case are then obtained in terms of  an analytic-geometric nondegeneracy condition, which generalize the result on the $C^{0}$ finite determination of  bifurcation problems given by P.B.Percell and P.N.Brown. \\}
 
 \end{minipage}
 \end{center}

 \section{{ Introduction}}
\ \
\indent \ \ Bifurcation phenomenon arise from a large number of nonlinear problems, the associated bifurcation problems are studied via a reduction ( such as the Lyapunov-Schmidt procedure ) to a finite dimensional local model $ G(u , \lambda):(\mathbf{R}^{n}\times \mathbf{R}^{l},0) \longrightarrow (\mathbf{R}^{p},0).$ This is viewed as a perturbation
$ G_{\lambda}(u)=G(u, ~\lambda)$ using parameters $\lambda\in  \mathbf{R}^{l}, $ of a germ $G_{0}(u)=G(u,0)$. There are a variety of notions of equivalence for studying such perturbations. The feature of interest in the present content is the variation of the set of zeros of $G_{\lambda}(0)$ with the parameter $\lambda.$

\indent A bifurcation problem is considered to be a family of maps
 $$ G(\cdot , \lambda):(\mathbf{R}^{n},0) \longrightarrow (\mathbf{R}^{p},0)$$
parameterized by $\lambda\in \mathbf{R}^{l}$ such that $G(0,0)=0$, or , more compactly , a map
$$ G(\cdot , \lambda):(\mathbf{R}^{n}\times \mathbf{R}^{l},0) \longrightarrow (\mathbf{R}^{p},0).$$
The set $G^{-1}(0)$ is called the bifurcation diagram. 
Roughly, two bifurcation problems are equivalent if their  bifurcation diagrams are locally homeomorphic  in a neighbourhood of the origin, and a very interesting problem is to determine what terms from the Taylor expansion at some points in neighbourhood of 0 may be omitted without changing the topological type determined by $G$ and the value of the bifurcation parameter $\lambda$. It is concerned with determinacy of bifurcation problems \\
\indent There is an extensive literature related to determination of bifurcation problems.
\indent $C^{\infty}$  theory of finite determinacy of bifurcation problems was systematically studied by M.Golubisky using singularity theory and group-theoretic techniques (ref [4], [5]). $C^{0}$ theory of finite determination of bifurcation problems is explored by Peter B.Percell and Peter N. Brown in [6]. They  have shown that $C^{0}$ finite determination of bifurcation diagrams follow from an analytic-geometric non-degeneracy condition which is modelled on a criterion of Kuo, rather than an algebraic condition of the type found in the  $C^{\infty}$ theory. In [14], Z. Jiangcheng, S.Fuwei, S.Ruixia and L.Guofu have discussed the d-determination of bifurcation problems with respect to $C^{0}$ contact equivalence from the weighted point of view, a criterion is given to judge the d-determination of bifurcation problems with respect to $C^{0}$ contact equivalence. Bucher, Marsden and Schecter ([2]) have also obtained a criterion for $C^{0}$ finite determination of bifurcation diagrams using a blowing-up construction and techniques from algebraic geometry. Above works on determination of bifurcation problems only deal with in the case that 
 perturbation $ G_{\lambda}(u)=G(u, ~\lambda)$ with parameter $\lambda$ of a germ $G_{0}(u)=G(u,0)$ has isolated bifurcation point in neighbourhood of the origin. \\
 \indent The case where a perturbation $ G_{\lambda}(u)=G(u, ~\lambda)$ with parameter $\lambda$ of a germ $G_{0}(u)=G(u,0)$ has non-isolated bifurcation points is more complicated. What we want to know is whether  $G$ is determined up to equivalence by finite coefficients in their Taylor expansion at every point that belongs to the subset of non-isolated bifurcation points of $ G_{\lambda}(u)=G(u, ~\lambda)$ with parameter $\lambda$ in neighbourhood of the origin.\\
\indent The bifurcation diagram of $G$ actually bifurcate at $(x~,~\lambda)$ it is  necessary 
$$rank~ [\bigtriangledown_{x}G(x~,~\lambda)]< p,$$ 
where $\bigtriangledown_{x}G$ denotes partial derivative of $G$ with respect to the variable $x\in\mathbf{R}^{n}.$ \\
\indent  Now if for a given closed set $\Sigma$ containing the origin in $(\mathbf{R}^{n},0),$
 there exists ~some~ neighbourhood ~$V$ of ~$0\in\mathbf{R}^{n}\times \mathbf{R}^{l}$ such that
\begin {eqnarray}
 rank~[\bigtriangledown_{x}G(x~,~\lambda)]=p ~~~(x~,~\lambda)\in V \setminus \Sigma \times \mathbf{R}^{l},
\end {eqnarray}
 then bifurcation points of $G$ only appear in set $\Sigma\times \mathbf{R}^{l}$ and the bifurcation diagram of $G$ has good  behaviour away from  $\Sigma\times \mathbf{R}^{l}.$ This leads to one to propose (1.1) as the basic nondegeneracy condition for bifurcation diagrams. However, since(1.1) alone is not always adequate as a criterion for finite determination, therefore our full nondegeneracy condition $(K_{\Sigma,\omega}^{r,\delta})_{\Gamma})$, which contain and refine condition  (1.1), will be stated in Section 2. It is a version of a condition of Kuo [15].\\
\indent  In this paper, we apply the idea of $r-\Sigma-C^{0}-$equivalence map jets relative to a given closed set $\Sigma$ in $(\mathbf{R}^{n},0),$ which is considered by Karim Bekka and Satoshi Koike [1] and   B. Osi$\acute{n}$ska-Ulrych, T. Rodak, G. Skalski in [13], to bifurcation theory, and introduce the $\Gamma-\Sigma-C^{0}-$bifurcation diagram and $\Gamma-\Sigma-C^{0}-$contact equivalence about $\Gamma-$equivariant bifurcation problems from the weighted point view and explore the finite determinacy of $\Gamma-$equivariant bifurcation problems with respect to above  equivalences. Some criteria about determination of $\Gamma-$equivariant bifurcation problems with respect to the $\Gamma-\Sigma-C^{0}-$bifurcation diagram and $\Gamma-\Sigma-C^{0}-$contact equivalence are given.\\

 \textbf{Theorem 1.1.}  {\sl Let $\Gamma$ be a compact Lie group acting orthogonally on space on space $\mathbf{R}^{n}\times \mathbf{R}^{l}$ and space $\mathbf{R}^{p}.~ 
f:(\mathbf{R}^{n}\times \mathbf{R}^{l},0) \longrightarrow (\mathbf{R}^{p},0)$ be a $\Gamma-$equiviant $C^{2}$ map-germ. Suppose $f$ satisfies the condition $(K_{\Sigma,~\omega}^{r,~\delta})_{\Gamma}$ and $\mid\omega\mid+\delta \geq 1,$ then $f$ is $ \Gamma-\Sigma-$ BD , or $ \Gamma-\Sigma-$ contact r-determined.}
\\

\indent The assumptions in the above theorem are natural and cannot be essentially improved.
 Our method used in proof of above theorem is concretely offering a controlled vector field. This vector field provides integration. But in order to obtain integration for giving a $\Gamma-\Sigma-C^{0}-$locally homeomorphic in a neighbourhood of the origin, it is necessary for a key tool to use. This tool is the existence and uniqueness of solution about differential equation in a singular Riemann metric on $\mathbf{R}^{n}\times \mathbf{R}^{l}$ from the weighted point view, which is stated and proved in section 2.\\
\indent When Lie group is a $\{e\}$, where $\{e\}$ is a unit element, and a singular Riemann metric transforms 
Euclidean metric, Theorem 1.1 implies \\ 
 
 \textbf{Theorem 1.2.}  {\sl Let $f:(\mathbf{R}^{n}\times \mathbf{R}^{l},0) \longrightarrow (\mathbf{R}^{p},0)$ be a $C^{r}$ map-germ and $p:(\mathbf{R}^{n}\times \mathbf{R}^{l},0) \longrightarrow (\mathbf{R}^{p},0)$ be a $c^{r}$ map-germ which satisfies:
 $$
 |p_{i}|=o(d(x, \Sigma)^{r}),~~~~~~|\frac{\partial p_{i}}{\partial x_{j}}|=o(d(x, \Sigma)^{r-1}),~~~i=1,\ldots,p;~~j=1,\ldots,n
 .$$
 If $f$ satisfies the condition $(K_{\Sigma}^{r,~\delta})$, then $f$ and $f+p$ is $ \Sigma-$BD r-determined.}\\
 
  \indent \textbf{Corollary 1.3.}  {\sl Suppose $~f:(\mathbf{R}^{n}\times \mathbf{R}^{l},0) \longrightarrow (\mathbf{R}^{p},0)$ be $C^{r+2}$ map-germ which satisfies the condition $(K_{\Sigma }^{r,\delta}))$ and $~g:(\mathbf{R}^{n}\times \mathbf{R}^{l},0) \longrightarrow (\mathbf{R}^{p},0)$ be $C^{r+2}$ map-germ which satisfies the following condition: there exists a neighbourhood U of $0$ in $\mathbf{R}^{n}\times \mathbf{R}^{l}$ such that, for every point $a\in (\Sigma\times \mathbf{R}^{l})\cap U,$ the jet $j^{\left( r+2\right) }f$ of Taylor formula of degree $\left( r+2\right)$ of $f$ at $a$ equals to the jet $j^{\left( r+2\right)}g$ of Taylor formula of $\left( r+2\right)$ of $g$ at $a,$ then $f$ and $g$ is $\Sigma-C^{0}$-BD equivalent and $\Sigma-C^{0}$-contact equivalent.}\\
  
  \indent Theorem 1.2 and Corollary 1.3 show to what terms from the Taylor expansion at every point that belongs to a closed subset $ \Sigma\times \mathbf{R}^{l}$ such that $(0,0)\in \Sigma\times \mathbf{R}^{l}$ may be omitted without changing the topological type determined by $f$ and the value of the bifurcation parameter $\lambda$.\\
  
 \indent Finally, we show that $C^{0}-$finite determination of  bifurcation problem, which is given by Peter B. Percell and Peter N. Brown In [6], is a corollary of Theorem 1.1.\\
 
  \textbf{Theorem 1.4.}([6], Theorem 3.1) {\sl Suppose $F\left( \mathbf{R}^{n}\times\mathbf{R}^{l},~0 \right)\rightarrow \left(  \mathbf{R}^{p},~0\right)$ is a $C^{1}$ map and $\nu=(\nu_{1},\cdots, ~\nu_{p})$ such that F is 
   ND$(\nu)$. Then F is contact $\mid\nu\mid-$ determined.}\\
   
 \indent  The rest of the paper is organized as follows. Section 2 contains necessary definitions, notation and technical preliminaries. Section 3 is devoted to proofs of Theorem 1.1, Theorem 1.2 and Corollary 1.3 respectively. In section 4, we give the proof of Theorem 1.4.
  
\section{{Preliminaries}}
\ \
\indent  Let $\Gamma$ be a compact Lie group. It acts linearly on $\mathbf{R}^{n}\times \mathbf{R}^{l}$ by definition
$$\gamma(x,~\lambda)=(\gamma x,~\lambda),~~~\forall \gamma \in \Gamma,~ \forall x\in \mathbf{R}^{n},~~\forall \lambda \in \mathbf{R}^{l}.$$
Meantime $\Gamma$ also acts linearly on $\mathbf{R}^{p},$  we say map-germ
$ G(\cdot , \lambda):(\mathbf{R}^{n}\times \mathbf{R}^{l},0) \longrightarrow (\mathbf{R}^{p},0)$ is $\Gamma-equiviant $ if
$$G(\gamma x, \lambda)=\gamma G(x,\lambda),~~\forall \gamma \in \Gamma,~ \forall x\in \mathbf{R}^{n},~~\forall \lambda \in \mathbf{R}^{l}$$
\indent let $ G,~F:(\mathbf{R}^{n}\times \mathbf{R}^{l},0) \longrightarrow (\mathbf{R}^{p},0)$ be $\Gamma-$equiviant continuous maps. We say $G$ and $F$ are $ \Gamma-[C^{r}]~ BD ~equivalent ~(BD~ for~ bifurcation ~diagram )$ if there are a neighbourhood $V\subset \mathbf{R}^{n}\times \mathbf{R}^{l}$ of the origin satisfying $\gamma\cdot V\subset V,~\forall \gamma \in \Gamma$
and a $\Gamma-$equiviant map persevering parameter level
$$\phi:~~(V, 0)\longrightarrow(\mathbf{R}^{n}\times \mathbf{R}^{l},0)$$ of the form
$$\phi(x,~\lambda)=(\phi_{1}(x,~\lambda),~\lambda),~~x\in \mathbf{R}^{n},~~\lambda\in \mathbf{R}^{l}$$
which is a homeomorphism [$C^{r}$ diffeomorphism]onto its image such that
$$\phi(G^{-1}(0)\bigcap V)=F^{-1}(0)\bigcap \phi (V).$$
We call $G$ and $F~are~\Gamma-[C^{r}] ~contact ~equivalent $ if there exist a neighbourhood $V\subset \mathbf{R}^{n}\times \mathbf{R}^{l}$ of the origin satisfying $\gamma\cdot V\subset V,~\forall \gamma \in \Gamma$
and a  $\Gamma-$equiviant map persevering parameter level
$$\phi:~~(V, 0)\longrightarrow(\mathbf{R}^{n}\times \mathbf{R}^{l},0)$$
which is  $ \Gamma-[C^{r}]$ BD equivalence between $G$ and $F$ and a continuous $[C^{r}]$ map
$$T:V\rightarrow GL(\mathbf{R}^{p})$$
satisfying
$$T(\gamma x,~\lambda)=\gamma\cdot T(x,~\lambda)\gamma^{-1},\gamma\in \Gamma$$
such that
$$G(u)=T(u)\cdot F(\phi(u)),~~~~u=( x,~\lambda).$$

\indent Let us fix a system of positive numbers $\omega=(\omega_{1},\ldots,\omega_{n},~\omega_{n+1},\cdots,\omega_{n+l})$, call $the~ weight ~of ~the ~variables~u_{i}~,~ i=1,~\ldots~,n+l.~\omega(u_{i})=\omega_{i}.$ and $(\mathbf{R}^{n}\times \mathbf{R}^{l},<\cdot>)$ be inner space. $\mathbf{e}_{1},\ldots \mathbf{e}_{n+l}$ is orthogonal basis of $(\mathbf{R}^{n}\times \mathbf{R}^{l},<\cdot>)$ .
For  $u=u_{1}\mathbf{e}_{1}+u_{2}\mathbf{e}_{2}+\ldots+u_{n+l}\mathbf{e}_{n+l}=x_{1}\mathbf{e}_{1}+\ldots+x_{n}\mathbf{e}_{n}+\lambda_{n+l}\mathbf{e}_{n+1}+\ldots+\lambda_{n+l}\mathbf{e}_{n+l}$ we may introduce the function
$$\rho=\rho(u)=\|u\|_{\omega}=(\sum^{n+l}_{i=1}u_{i}^{2q_{i}})^{\frac{1}{2q}}
=(\sum^{n}_{i=1}x_{i}^{2q_{i}}+\sum^{n+l}_{i=n+1}\lambda_{i}^{2q_{i}})^{\frac{1}{2q}},$$
where $q_{i}=\frac{q}{\omega_{i}},~~1\leq i\leq n+l $ and $q=\omega_{1}\omega_{2}\cdots\omega_{n+l}.$ \\

\textbf{ Remark 1.}([13]) $\rho (u)^{2q}$ satisfies a Lojasiewicz condition $\rho (u)^{2q}\geq c \parallel u\parallel^{2\alpha}$ for some constants $c$ and $\alpha.$\\

 \textbf{Definition 2.1}([3]) Using this $\rho$ we may introduce the $ singular~ Riemannian~ metric~$ on $(\mathbf{R}^{n}\times \mathbf{R}^{l},<\cdot>)$, namely the Riemannian metric on  $(\mathbf{R}^{n}\times \mathbf{R}^{l}\setminus\{0\})$ defined by the following bilinear form:
 $$\langle\frac{\partial}{\partial u_{i}},~\frac{\partial}{\partial u_{i}}\rangle=\rho^{-2\omega_{i}},~~~\langle\frac{\partial}{\partial u_{i}},~\frac{\partial}{\partial u_{j}}\rangle=0,~~1\leq i,~j\leq n+l,~~i\neq j.$$
 $$
 \langle du_{i_{1}}\bigwedge\cdots\bigwedge du_{i_{k}},~ du_{i_{1}}\bigwedge\cdots\bigwedge du_{i_{k}}\rangle=\|u\|_{\omega}^{2(\omega_{i_{1}}+\cdots+\omega_{i_{k}})}=\rho^{2(\omega_{i_{1}}+\cdots+\omega_{i_{k}})}.
 $$
 We denote by $\nabla_{\omega},~\|\cdot\|_{\omega},$ the corresponding gradient and norm associated with this Riemannian metric. For a function-germ
$ G,F:(\mathbf{R}^{n}\times \mathbf{R}^{l},0) \longrightarrow (\mathbf{R},0),$
$$
\nabla_{\omega}G =\sum^{n+l}_{i=1}\rho^{\omega_{i}}\frac{\partial G}{\partial u_{i}}\rho^{\omega_{i}}\frac{\partial }{\partial u_{i}}
,~~~~
\|\nabla_{\omega}G\|_{\omega}^{2}=\sum^{n+l}_{i=1}(\rho^{\omega_{i}}\frac{\partial G}{\partial u_{i}})^{2},
$$
$$\left\langle \nabla_{\omega}G,~\nabla_{\omega}F \right\rangle =
\sum^{n+l}_{i=1}\rho^{2\omega_{i}}\frac{\partial G}{\partial u_{i}}\frac{\partial F }{\partial u_{i}}.
$$
If a map-germ
$ G(\cdot , \lambda):(\mathbf{R}^{n}\times \mathbf{R}^{l},0) \longrightarrow (\mathbf{R}^{p},0),$ then the gradient of component $g_{i}$ of $G$ with respect to $x$ is

$$
\nabla_ {\omega,x}g_{i} =\sum^{n}_{i=1}\rho^{\omega_{j}}\frac{\partial g_{i}}{\partial x_{j}}\rho^{\omega_{j}}\frac{\partial }{\partial x_{j}}
,~~~~~
\|\nabla_{\omega,x}g_{i}\|_{\omega}^{2}=\sum^{n}_{j=1}(\rho^{\omega_{j}}\frac{\partial g_{i} }{\partial x_{j}})^{2},
$$

We still denote by $\mathbf{R}^{n}\times \mathbf{R}^{l}$ the inner linearly space $\mathbf{R}^{n}\times \mathbf{R}^{l}$ with this singular metric.\\
\indent$ The~ weighted ~horn~neighbourhood $ of degree $d$ and width $c>0$ of a variety $G^{-1}(0) $ is by definition
$$H_{d}(G,~c)=\{u\in \mathbf{R}^{n}\times \mathbf{R}^{l}|~~ \|G(u)\|_{(\mathbf{R}^{p}, <\cdot>)}\leq c\rho^{d}\}.$$

\indent Now let $\Sigma$ be a germ of a closed subset of $\mathbf{R}^{n}$ such that $0\in \Sigma.$ Then we denote by $\mathcal{R}_{\Sigma\times \mathbf{R}^{n}}^{fix}$ the group of germs of homeomorphisms   $\phi:~(\mathbf{R}^{n}\times \mathbf{R}^{l},~0)\rightarrow (\mathbf{R}^{n}\times \mathbf{R}^{l},~0)$ at $0\in \mathbf{R}^{n}\times \mathbf{R}^{l}$ which fixes $\Sigma\times \mathbf{R}^{l}$ namely
$$\phi(x,\lambda)=((\phi_{1}(x,\lambda),\lambda)=(x,\lambda),~~~\forall x\in \Sigma,~~\lambda \in\mathbf{R}^{l}$$

\indent We consider the following equivalence relation: \\

  \textbf{Definition 2.2}([6]) Let $\Gamma$ be a compact Lie group, $~ G,~F:(\mathbf{R}^{n}\times \mathbf{R}^{l},0) \longrightarrow (\mathbf{R}^{p},0)$ be $\Gamma-$equiviant continuous maps. We say that $G$ and $F$ are $ \Gamma-\Sigma-[C^{r}]~ BD~ equivalent~ \\
 (BD~ for~ bifurcation~ diagram )$ if there are a neighborhood $V\subset \mathbf{R}^{n}\times \mathbf{R}^{l}$ of the origin satisfying $\gamma\cdot V\subset V,~\forall \gamma \in \Gamma$
and a homeomorphism $\phi \in\mathcal{R}_{\Sigma\times \mathbf{R}^{n}}^{fix}$
 [$C^{r}$ diffeomorphism]onto its image such that
$$\phi(G^{-1}(0)\bigcap V)=F^{-1}(0)\bigcap \phi (V)$$
and fixes $G^{-1}(0)\bigcap \Sigma.$\\
\indent We call that $G$ and $F $ are $\Gamma-\Sigma-[C^{r}]~ contact ~equivalent$ if there exist a neighbourhood $V\subset \mathbf{R}^{n}\times \mathbf{R}^{l}$ of the origin satisfying $\gamma\cdot V\subset V,~\forall \gamma \in \Gamma$
and a homeomorphism $\phi \in\mathcal{R}_{\Sigma\times \mathbf{R}^{n}}^{fix}$[$C^{r}$ diffeomorphism]
which is  $ \Gamma-\Sigma-[C^{r}]$ BD equivalence between $G$ and $F$ and a continuous $\Gamma-[C^{r}]$ equivariant map
$$T:V\rightarrow GL(\mathbf{R}^{p}),$$
namely $$ T(\phi_{1}(x,\lambda),\lambda)=\gamma T(x,\lambda)\gamma^{-1},~~~\forall \gamma\in \Gamma,$$
such that
$$G(u)=T(u)\cdot F(\phi(u)),~~~~u=( x,~\lambda),~~u=(x,\lambda)\in V.$$

\indent Now we introduce the finite determination concepts.\\

\indent  \textbf{Definition 2.3}\ \ Let $\Gamma$ be a compact Lie group and  $f:(\mathbf{R}^{n}\times \mathbf{R}^{l},0) \longrightarrow (\mathbf{R}^{p},0)$ be $\Gamma-$equiviant continuous map. We say that $f$ is $ \Gamma-\Sigma-[C^{r}]~ BD~ ,~ or ~ \Gamma-\Sigma-[C^{r}]~ contact~ d-determined $ if $f$ and $f+p$ are  $ \Gamma-\Sigma-[C^{r}]$ BD , or $ \Gamma-\Sigma-[C^{r}]$ contact equivalent for every $C^{2}[C^{max(r,~2)}]~\Gamma-$ equivariant perturbation
$$p: (\mathbf{R}^{n}\times \mathbf{R}^{l})\rightarrow (\mathbf{R}^{p},0)$$
such that, with $u=(x,\lambda)\in\mathbf{R}^{n}\times \mathbf{R}^{l}$ and
\begin{eqnarray}
|p_{i}|=o(d_{\omega}(x, \Sigma)^{d+\vert\omega\vert}),~~~~~~|\frac{\partial p_{i}}{\partial x_{j}}|=o(d_{\omega}(x, \Sigma)^{d}),~~~i=1,\ldots,p;~~j=1,\ldots,n,
\end{eqnarray}
where $$
d_{\omega}(x, \Sigma)=inf_{y\in \Sigma}\left\lbrace \parallel x-y\parallel_{\omega} \right\rbrace 
$$ 
and 
$$\vert\omega\vert=max\{\vert\omega_{i}\vert:~i=1,\cdots,~n+l\}.$$\\

 \indent Let $\Sigma$ be a germ of closed set of $\mathbf{R}^{n}$ such that $0\in \Sigma.$
 Let $G$ be a germ of $C^{1}$ vector field on $\mathbf{R}^{n}\times \mathbf{R}^{l}\setminus \Sigma\times \mathbf{R}^{l} $ which satisfies the relative Lipschitz condition: $ \parallel G(x, t)\parallel \leq C d(x,~\Sigma).$ \\
 
 \indent For a fixed vector $v\in \lbrace 0\rbrace \times \mathbf{R}^{l},$ we define
  $$
 X(u,~t)= \left\{
 \begin {array}{cc}
 G(x,~t)+v,~~~if~ x\in \mathbf{R}^{n}\setminus \Sigma.
 \\[6pt]
v,~~~~~if~ x\in \Sigma\\[6pt]
 \end {array}
        \right.
 $$
\indent Then  \\

  \indent  \textbf{Lemma 2.4}([1], Proposition 2,15 ) {\sl For~  G(x,~t)~ satisfying~ the~ preceding ~conditions~$,~ $ X(u,~t)~ is~ locally ~integrable~ in ~the sense ~that ~there~ are ~a ~neighbourhood ~W ~of ~$(x_{0},~t_{0})$~ in~ ~$ \mathbf{R}^{n}\times \mathbf{R}^{l},\delta> 0,$~ and ~a ~family~ of~ homeomorphisms ~$ \Phi_{s}~$
defined ~on ~W~ for~ $ s<\delta$
so~ that~$ \Phi_{0}=id $~and~ for~ $(x,~t,~s)\in W\times (-\delta,~\delta),$ 
$$
\frac{\partial \Phi_{s}} {\partial s}=X\circ\Phi_{s} .
$$}

  \indent  \textbf{Lemma 2.5}([1], Lemma 2.16; [13],Lemma 2.13)  {\sl Let $U$ be an open subset of $ \mathbf{R}^{n}\setminus \Sigma ,$
  let $0\in (a,~b)$ and let $G:~U\times (a,~b)\rightarrow  \mathbf{R}^{n}$ be a continuous mapping which satisfies
   $$ \parallel G(x, t)\parallel \leq C d(x,~\Sigma)$$ 
   for some $C>0$ and $(x,~t)\in U\times (a,~b).$ Let $\varphi (\alpha,~\beta)\rightarrow U$ be an integral solution of the system of differential equations $y^{'}=G(y,~t)$
with the initial condition $\varphi(0)=x_{0}$ where $x_{0}\in U$ and $0\in (\alpha,~\beta)\subset (a,~b).$
Then we have 
 \begin {eqnarray} 
 d (x_{0},~\Sigma)e^{-C\mid t \mid}\leq  d (\varphi (t),~\Sigma)\leq d (x_{0},~\Sigma)e^{C\mid t \mid}
 \end {eqnarray}
for $t\in (\alpha,~\beta). $}\\

\indent Next we provide the following key proposition which is used in the proof of Theorem 1.1.\\

 \indent  \textbf{Proposition 2.6} (Key Proposition){\sl For~  G(x,~t)~ satisfying~
 $$ \parallel G(x, t)\parallel_{\omega} \leq C d_{\omega}(x,~\Sigma),$$ 
 ~ $ X(u,~t)$~ is~ locally integrable~ in ~the sense ~that ~there~ are ~a ~neighbourhood $~W$ ~of ~$(x_{0},~t_{0})~ \in~ ~ \mathbf{R}^{n}\times \mathbf{R}^{l},\delta> 0,$~ and ~a ~family~ of~ homeomorphisms ~$ \Phi_{s}$
defined ~on $~W$~ for~ $ s<\delta$
so~ that~ $\Phi_{0}=id$ ~and~ for~ $(x,~t,~s)\in W\times (-\delta,~\delta),$ 
$$
\frac{\partial \Phi_{s}} {\partial s}=X\circ\Phi_{s} .
$$
 \\
 \indent  \textbf{Lemma 2.7} {\sl Let $U$ be an open subset of $ \mathbf{R}^{n}\setminus \Sigma ,$
  let $0\in (a,~b)$ and let $G:~U\times (a,~b)\rightarrow  \mathbf{R}^{n}$ be a continuous mapping which satisfies
   \begin {eqnarray}
    \parallel G(x, t)\parallel_{\omega} \leq C d_{\omega}(x,~\Sigma)
     \end {eqnarray} 
   for some $C>0$ and $(x,~t)\in U\times (a,~b).$ Let $\varphi (\alpha,~\beta)\rightarrow U$ be an integral solution of the system of differential equations $y^{'}=G(y,~t)$
with the initial condition $\varphi(0)=x_{0}$ where $x_{0}\in U$ and $0\in (\alpha,~\beta)\subset (a,~b).$
Then we have 
 \begin {eqnarray*} 
d_{\omega} (\varphi (0),~\Sigma) e^{-CL\mid t \mid}<  d_{\omega} (\varphi (t),~\Sigma)\leq d_{\omega} (\varphi (0),~\Sigma) e^{CL\mid t \mid},
 \end {eqnarray*}
 for $t\in (\alpha,~\beta). $} and some positive number $L$.}\\
 
\indent  {\sl Proof.} Owing to $ \varphi( t)\in U,$ for $y\in \Sigma$ and $t\in (\alpha,~\beta),$ let $\parallel \varphi_{y_{•}}( t)\parallel_{\omega}=\parallel \varphi(t)-y\parallel_{\omega}>0$ and  $\varphi_{y_{•}}( t)= \left(\varphi_{y,1}( t),~\cdots,~\varphi_{y,~n}( t) \right) $\\
\indent We can get the function $\kappa(\varphi_{y}( t))=\frac{1}{2\omega}ln\parallel \varphi_{y}( t)\parallel_{\omega}^{2\omega} =\frac{1}{2\omega}ln\rho \left( \varphi_{y}( t)\right)^{2\omega}$ for $t\in (\alpha,~\beta)$
and 
\begin {eqnarray*}
\frac{d\kappa(\varphi_{y}( t))}{dt}&=&\frac{1}{2\omega}\frac{(2\omega)\parallel \varphi_{y}( t)\parallel_{\omega}^{2\omega-1}}{\parallel \varphi_{y}( t)\parallel_{\omega}^{2\omega}}\cdot\frac{d\parallel \varphi_{y}( t)\parallel_{\omega}}{dt}\\
&=&\frac{1}{\parallel \varphi_{y}( t)\parallel_{\omega}}\cdot\frac{d\parallel \varphi_{y}( t)\parallel_{\omega}}{dt}\\
&=&\frac{1}{\parallel \varphi_{y}( t)\parallel_{\omega}}\cdot\sum_{i=1}^{n}\frac{\partial\rho(x)}{\partial x_{i}}\cdot
\frac{\varphi_{y,i}( t)}{dt}\\
&=&\frac{1}{\parallel \varphi_{y}( t)\parallel_{\omega}}\cdot\sum_{i=1}^{n}\frac{\partial\rho(x)}{\partial x_{i}}\cdot G_{i}\left( \varphi_{y}( t),~t\right) ~~~for ~t\in (\alpha,~\beta).
\end {eqnarray*}

\indent Since 
\begin {eqnarray*}
\frac{\partial\rho(x)}{\partial x_{i}}&=&\frac{1}{2\omega}(\sum^{n}_{i=1}x_{i}^{\frac{2\omega}{\omega_{i}}})^{\frac{1}{2\omega}-1}\cdot\frac{2\omega}{\omega_{i}}\cdot x_{i}^{\frac{2\omega}{\omega_{i}}-1}\\
&=&\frac{1}{\omega_{i}}\frac{\rho\cdot x_{i}^{\frac{2\omega}{\omega_{i}}}}{ (\sum^{n}_{i=1}x_{i}^{\frac{2\omega}{\omega_{i}}})}\cdot\frac{1}{x_{i}},
\end {eqnarray*}
$$
\frac{\partial\rho(x)}{\partial x_{i}}\cdot\rho^{\omega_{i}-1}=\frac{1}{\omega_{i}}\frac{ x_{i}^{\frac{2\omega}{\omega_{i}}}}{ (\sum^{n}_{i=1}x_{i}^{\frac{2\omega}{\omega_{i}}})}\cdot\frac{(\sum^{n}_{i=1}x_{i}^{\frac{2\omega}{\omega_{i}}})^{\frac{\omega_{i}}{2\omega}}}{x_{i}}
$$
\indent Now let us observe that
\begin {eqnarray}
\mid\frac{\partial\rho(x)}{\partial x_{i}}\mid\leq \frac{L}{n}\rho^{1-\omega_{i}}~~for ~some ~L>0~~(also~by~ ([7])).
\end {eqnarray}
and
\begin {eqnarray*}
\mid~G_{i}\left( \varphi_{y}( t),~t\right)\mid=\left[\mid~G_{i}\left( \varphi_{y}( t),~t\right)\mid^{\frac{1}{•\omega_{i}}} \right]^{\omega_{i}} & \leq & \left[ \parallel ~G\left( \varphi_{y}( t),~t\right)\parallel_{\omega}\right] ^{\omega_{i}} \\
&\leq & C \left[ \parallel\varphi_{y}( t)\parallel_{\omega}\right] ^{\omega_{i}}~~(by~(2.2)).
\end {eqnarray*}
I.e.
\begin {eqnarray}
\mid~G_{i}\left( \varphi_{y}( t),~t\right)\mid \leq  C \left[ \parallel\varphi_{y}( t)\parallel_{\omega}\right] ^{\omega_{i}}.
\end {eqnarray}

\indent From the mean value theorem, for every $t\in (0,~\beta),$ there exists $\theta\in (0,~\beta)$
such that 
\begin {eqnarray*}
\mid\kappa( t)-\kappa( 0)\mid &\leq &\mid\frac{d\kappa(\varphi_{y}( \theta))}{dt}\mid t
=\frac{1}{\parallel \varphi_{y}( t)\parallel_{\omega}}\cdot\mid\sum_{i=1}^{n}\frac{\partial\rho(x)}{\partial x_{i}}\cdot G_{i}\left( \varphi_{y}( t),~t\right)\mid t\\
&\leq &\frac{1}{\parallel \varphi_{y}( t)\parallel_{\omega}}\cdot\sum_{i=1}^{n}\mid \frac{\partial\rho(x)}{\partial x_{i}}\mid\cdot\mid~G_{i}\left( \varphi_{y}( t),~t\right)\mid  t\\
&\leq &\frac{1}{\parallel \varphi_{y}( t)\parallel_{\omega}}\cdot\sum_{i=1}^{n}\frac{L}{n}\rho^{1-\omega_{i}}\cdot C \left[ \parallel\varphi_{y}( t)\parallel_{\omega}\right] ^{\omega_{i}} t\\
&=&  C\sum_{i=1}^{n}\frac{1}{\parallel \varphi_{y}( t)\parallel_{\omega}^{1-\omega_{i}}}\frac{L}{n}\rho^{1-\omega_{i}}\cdot  t=C \cdot L\cdot t,
\end {eqnarray*}
where $\rho=\parallel \varphi_{y}( t)\parallel_{\omega}.$

\indent Hence for every $t\in (0,~\beta),$
$$
\kappa( 0)-C L t\leq\kappa( t)\leq\kappa( 0)+C L t.
$$
The above inequalities hold also for $t=0,$ i.e.
\begin {eqnarray}
\parallel \varphi_{y}( 0)\parallel_{\omega}e^{-C L\mid t\mid}\leq\parallel \varphi_{y}( t)\parallel_{\omega}
\leq\parallel \varphi_{y}( 0)\parallel_{\omega}e^{C L\mid t\mid}.
\end {eqnarray}
\indent Because $ d_{\omega} (\varphi (t),~\Sigma)=inf_{y\in \Sigma}{\parallel \varphi_{y}( t)\parallel_{\omega}}, $ and (2.4), we obtain
\begin {eqnarray} 
 d_{\omega} (\varphi (0),~\Sigma) e^{-CL\mid t \mid}<  d_{\omega} (\varphi (t),~\Sigma)\leq  d_{\omega} (\varphi (0),~\Sigma) e^{CL\mid t \mid}.
 \end {eqnarray}\\
 
{\sl Proof of Proposition 2.6.} the proof will be essentially the same as that of Proposition 2.15 of [1]  using Lemma 2.7.\\

\indent  \textbf{Lemma 2.8.} ([1], Lemma 2.4 ) {\sl Let $\Sigma$ be a germ at $0\in \mathbf{R}^{n} $ of a closed subset, and $f:~ (\mathbf{R}^{n},0) \longrightarrow (\mathbf{R}^{p},0)$ be a $C^{k}$ map-germ, $k\geq 1,$ which satisfies the following condition:\\
\indent there exists a neighbourhood U of $0$ in $\mathbf{R}^{n}$ such that, for every point $a\in \Sigma\cap U,$ the k-jet $j^{k}f(a)$ of Taylor formula of degree k of $f$ at $a$ is $0.$ Then $\parallel f(x)\parallel=o(d(k,~\Sigma)^{k}).$}\\
 
\indent  Let
 map $ f:(\mathbf{R}^{n}\times \mathbf{R}^{l},0) \longrightarrow (\mathbf{R}^{p},0)$ is a $C^{r}$ map,
 we consider vectors
   $$N(f,j,u)=\nabla_{x}f_{j}(u,t)-q_{j}(u,t)~~1\leq j\leq p,$$ 
  where $q_{j}(u,t)$ is the projection of $\nabla_{x}f_{j}(u,t)$ to the subspace $V_{u}^{j}$ spanned by $\left\lbrace \nabla_{x}f_{i}\right\rbrace _{j\neq i}.$ \\
   {\sl The Kuo pseud-distance } $d_{\omega, x}\nabla f$ is defined by 
  $$
 d_{\omega,x}\nabla f=min_{ 1 \leq i \leq p}\left\lbrace \|N(f,i,x)\|_{\omega}\right\rbrace .  
 $$
 
\indent  \textbf{Definition 2.9.} (the relative Kuo condition $(K_{\Sigma}^{r,\delta})$)   Let
 map $ f:(\mathbf{R}^{n}\times \mathbf{R}^{l},0) \longrightarrow (\mathbf{R}^{p},0)$ is a $C^{2}$ map,
 the map $f$ satisfies {\sl the relative Kuo condition} $(K_{\Sigma}^{r,\delta})$ if there is a strictly positive number $C,~\delta, \alpha$ and $\overline{w}$ such that
$$
d_{\omega,x}\nabla f\geq C\cdot d(x,~\Sigma)^{r-\delta}
$$
holds on $ u=(x,~\lambda)\in   H_{r}(f,~\overline{w})\bigcap\{\|u\|<\alpha\} ,$
where
$$
H_{r}^{\Sigma}(f,~\overline{w})=\{u\in \mathbf{R}^{n}\times \mathbf{R}^{l}:\|f(u)\|\leq \overline{w}d(x,~\Sigma)^{r} \}.
$$
\indent  \textbf{Definition 2.10.}(the relative Kuo condition$(K_{\Sigma,\omega}^{r,\delta})_{\Gamma})$)  Let $\Gamma$ be a compact Lie group acting on space on space $\mathbf{R}^{n}\times \mathbf{R}^{l}$ and space $\mathbf{R}^{p}.$ and $\Sigma$ be a germ of closed set of $\mathbf{R}^{n}$ such that $0\in \Sigma$ and  $\gamma\cdot\Sigma\subset \Sigma,~\forall \gamma \in \Gamma$. A $C^{2}~\Gamma-$equivariant map $ f:(\mathbf{R}^{n}\times \mathbf{R}^{l},0) \longrightarrow (\mathbf{R}^{p},0)$ is  said to satisfy {\sl the relative Kuo condition} $(K_{\Sigma}^{r,\delta})_{\Gamma}$ if there is a strictly positive number $C,~\delta, \alpha$ and $\overline{w}$ such that
$$
d_{\omega,x}\nabla f\geq C\cdot d_{\omega}(x,~\Sigma)^{r-\delta}
$$
holds on $ u=(x,~\lambda)\in   H_{\omega,r}^{\Sigma}(f,~\overline{w})\bigcap\{\|u\|<\alpha\} ,$
where
$$
H_{\omega,r}^{\Sigma}(f,~\overline{w})=\{u\in \mathbf{R}^{n}\times \mathbf{R}^{l}:\|f(u)\|\leq \overline{w}d_{\omega}(x,~\Sigma)^{r} \}.
$$\\
\textbf{ Remark 2.}Since $\rho (u)^{2q}\geq c \parallel u\parallel^{2\alpha}$ for some constants $c$ and $\alpha$ by  Remark 1, then
$$
H_{\frac{ \alpha r}{q}}^{\Sigma}(f,~\overline{w}c^{\frac{ r}{2q}})\subset H_{\omega,r}^{\Sigma}(f,~\overline{w}).
$$
In fact, if $ u\in H_{\frac{ \alpha r}{q}}^{\Sigma}(f,~\overline{w}c^{\frac{ r}{2q}}),$
$$
\|f(u)\|\leq \overline{w}c^{\frac{ r}{2q}}d(x,~\Sigma)^{\frac{ \alpha r}{q}}=\overline{w}\left(cd(x,~\Sigma)^{2\alpha} \right)^{\frac{ r}{2q}} \leq \overline{w}d_{\omega}(x,~\Sigma)^{r}.
$$

\section{{The determinacy of $\Gamma-$equivariant bifurcation problems with respect to $ \Gamma-\Sigma-$ BD and $ \Gamma-\Sigma-$contact equivalence from the weighted point view}}
\ \

 \indent  In order to prove Theorem 1.1, we need to following Lemma.\\

\ \textbf{Lemma 3.1. } 
  Let $\Gamma$ be a compact Lie group acting orthogonally on space on space $(\mathbf{R}^{n}\times \mathbf{R}^{l},~\langle\cdot\rangle)$ and
 space $(\mathbf{R}^{p},~\langle\cdot\rangle),$ then,\\
  \indent (1) the representation of $\Gamma$ on the inner space $(\mathbf{R}^{n}\times \mathbf{R}^{l},\langle \cdot \rangle)$ is orthogonal, i.e.$\langle \gamma x,~\gamma y\rangle_{\Gamma}=\langle x,~y\rangle_{\Gamma }
.$\\
 \indent (2) $\rho (u)=\parallel u \parallel_{\omega}=\parallel \gamma u \parallel_{\omega}=\rho (\gamma u),~~~\forall \gamma \in \Gamma$. \\
  \indent (3) if $~f:(\mathbf{R}^{n}\times \mathbf{R}^{l},0) \longrightarrow (\mathbf{R}^{p},0)$ is $\Gamma-$equiviant $C^{r}$ map-germ and $\Sigma$ be a germ of closed set of $\mathbf{R}^{n}$ such that $0\in \Sigma$ and  $\gamma\cdot\Sigma\subset \Sigma,~\forall \gamma \in \Gamma$, then, for $\gamma\in \Gamma,~\gamma\cdot H_{r}^{\Sigma}(f,~\overline{w})=H_{r}^{\Sigma}(f,~\overline{w})$ and $\{\|\gamma u\|<\alpha\}=\{\|u\|<\alpha\}$ and $~\gamma\cdot H_{\omega,r}^{\Sigma}(f,~\overline{w})=H_{\omega,r}^{\Sigma}(f,~\overline{w}).$ I.e. $ H_{r}^{\Sigma}(f,~\overline{w})$ and $ H_{\omega,r}^{\Sigma}(f,~\overline{w})$ are $\Gamma-$ invariant set.\\
   
 \indent {\sl Proof.} (1) is obvious.\\
  \indent (2) In fact, for $\gamma\in \Gamma,$ if the representation of $\gamma$ is orthogonal matrix $A(\gamma)=(a_{ij})$, then 
 $$
 \gamma\cdot u=A(\gamma)\cdot u =(v_{1},v_{2},\cdots,v_{n+l} )=(\Sigma_{1\leq i\leq n+l}a_{1i}u_{i},\cdots,\Sigma_{1\leq i\leq n+l}a_{n+l ~i}u_{i})
 $$
By Definition 2.1, 
 $$
 \langle dv_{i_{1}}\bigwedge\cdots\bigwedge dv_{i_{k}},~ dv_{i_{1}}\bigwedge\cdots\bigwedge dv_{i_{k}}\rangle=\|\gamma\cdot u\|_{\omega}^{2(\omega_{i_{1}}+\cdots+\omega_{i_{k}})}=\rho(v)^{2(\omega_{i_{1}}+\cdots+\omega_{i_{k}})}.
 $$
So
$$
 \langle dv_{1}\bigwedge\cdots\bigwedge dv_{n+l},~ dv_{1}\bigwedge\cdots\bigwedge dv_{n+l}\rangle=\|\gamma\cdot u\|_{\omega}^{2\left( \omega_{1}+\cdots+\omega_{n+l}\right) }=\rho(v)^{2(\omega_{1}+\cdots+\omega_{n+l})}.
 $$
Moreover 
 \begin {eqnarray*}
 &~& \langle  dv_{1}\bigwedge\cdots\bigwedge dv_{n+l},~ dv_{1}\bigwedge\cdots\bigwedge dv_{n+l}\rangle=\\
  & = & \langle  \mid A(\gamma)\mid d u_{1}\bigwedge\cdots\bigwedge du_{n+l},~~ \mid A(\gamma)\mid d u_{1}\bigwedge\cdots\bigwedge du_{n+l}\rangle \\
  & = & \langle d u_{1}\bigwedge\cdots\bigwedge d u_{n+l},~ d u_{1}\bigwedge\cdots\bigwedge d u_{n+l}\rangle \\
 & = & \| u\|_{\omega}^{2(\omega_{1}+\cdots+\omega_{n+l})}=\rho(u)^{2(\omega_{1}+\cdots+\omega_{n+l})}
 \end {eqnarray*}
 So
  $$
  \rho(u)^{2(\omega_{1}+\cdots+\omega_{n+l})}=\rho(\gamma \cdot u)^{2(\omega_{1}+\cdots+\omega_{n+l})}
 $$
and $\rho(u)=\rho(\gamma u).$\\
 \indent (3) Since $f$ is $\Gamma-$equivariant and $\Gamma$ is orthogonally act on $\mathbf{R}^{n}\times \mathbf{R}^{l}$ and $\mathbf{R}^{p},$ then $$
 \parallel f(\gamma u)\parallel=\parallel\gamma f( u)\parallel=\parallel f( u)\parallel 
,$$
$$
d(\gamma x,~\Sigma)=d(\gamma x,~\gamma\Sigma)=d( x,~\Sigma).
$$
and
$$
d_{\omega}(\gamma x,~\Sigma)=d_{\omega}(\gamma x,~\gamma\Sigma)=d_{\omega}( x,~\Sigma).
$$
$~\parallel f( u)\parallel \leq \overline{w}d( x,~\Sigma)^{r}$ if and only if $\parallel f(\gamma\cdot u)\parallel \leq \overline{w}d(\gamma x,~\Sigma)^{r}=\overline{w}d( x,~\Sigma)^{r}$ and
$~\parallel f( u)\parallel \leq \overline{w}d_{\omega}( x,~\Sigma)^{r}$ if and only if $\parallel f(\gamma\cdot u)\parallel \leq \overline{w}d_{\omega}(\gamma x,~\Sigma)^{r}=\overline{w}d_{\omega}( x,~\Sigma)^{r}.$\\
\indent Therefore , for $\gamma\in \Gamma,~\gamma\cdot H_{r}^{\Sigma}(f,~\overline{w})\subset H_{r}^{\Sigma}(f,~\overline{w}) $ and $\gamma\cdot H_{\omega,r}^{\Sigma}(f,~\overline{w})\subset H_{\omega,r}^{\Sigma}(f,~\overline{w}).$ Again since $u=\gamma\cdot\left( \gamma^{-1}u\right),~  H_{r}^{\Sigma}(f,~\overline{w})\subset ~\gamma\cdot H_{r}^{\Sigma}(f,~\overline{w}) $ and $   H_{\omega,r}^{\Sigma}(f,~\overline{w})\subset ~\gamma\cdot H_{\omega,r}^{\Sigma}(f,~\overline{w}) .$ I.e. $ H_{r}^{\Sigma}(f,~\overline{w})$ and $ H_{\omega,r}^{\Sigma}(f,~\overline{w})$ are $\Gamma-$ invariant set.\\

\textbf{Lemma 3.2.}{ \sl Let $\Gamma$ be a compact Lie group acting orthogonally on space on space $(\mathbf{R}^{n}\times \mathbf{R}^{l},~\langle\cdot\rangle)$ and
 space $(\mathbf{R}^{p},~\langle\cdot\rangle)$. Suppose $~f,p:(\mathbf{R}^{n}\times \mathbf{R}^{l},0) \longrightarrow (\mathbf{R}^{p},0)$ is $\Gamma-$equiviant $C^{r}$ map-germ such that $p(u)$ satisfies (2.1).
  Define
  $$F(u,~t)=f(u)+tp(u),~~u=(x,~\lambda)\in \mathbf{R}^{n}\times \mathbf{R}^{l},t\in I=[0,~1].$$

  If there is a strictly positive number $C,~\delta, ~\alpha$ and $\overline{w}$ such that  $f$ satisfies the condition $(K_{\Sigma}^{r, \delta})_{\Gamma}$ as above, then when
 $$
   (u,~t)=(x,~\lambda,~t)\in H_{\omega,r}^{\Sigma}(f,~\overline{w})\bigcap\{\|u\|<\alpha\}\times I,
   $$there exists some positive constant $C^{\prime}$ such that
$$ d_{\omega,x}\nabla F=d_{\omega,x}(\nabla_{\omega,x}F_{1},\ldots,~\nabla_{\omega,x}F_{p})\geq C^{\prime} d_{\omega}(x,~\Sigma)^{r-\delta}$$ holds.}\\

{\sl Proof.} First we prove 
\begin{eqnarray}
 \| \nabla_{\omega,x}F_{i}(u,t)-\nabla_{\omega,x}f_{i}(u)\|_{\omega}\leq o(d_{\omega}(x, \Sigma)^{r})
\end{eqnarray}
In fact, by the definition of the singular Riemannian metric,
 \begin{eqnarray*}
    \| \nabla_{\omega,x}F_{i}(u,t)-\nabla_{\omega,x}f_{i}(u)\|_{\omega} &=& \| \nabla_{\omega,x}tp_{i}(u)\|_{\omega}\\ 
      & =&\left(t^{2}\sum^{n}_{j=1}(\rho^{\omega_{j}}\frac{\partial p_{i}}{\partial x_{j}})^{2} \right)^{\frac{1}{2}} \\
      & \leq & \mid t\mid(\sum^{n}_{j=1}\rho^{\omega_{j}}\mid\frac{\partial p_{i}}{\partial x_{j}}\mid)\\
     & =&\mid t\mid(\sum^{n}_{j=1}\rho^{\omega_{j}}o(d_{\omega}(x,~\Sigma)^{r}) \\ 
     & =& o(d_{\omega}(x,~\Sigma)^{r}) .
     \end{eqnarray*}
 So (3.7)holds.\\
 \indent  Now for $(u,~t)=(x,~\lambda,~t)\in \{\parallel u\parallel<\alpha\}\backslash \Sigma\times \mathbf{R}^{l}\times I,~ $ the vectors $\{\nabla_{\omega,x}f_{1},\ldots, \nabla_{\omega,x}f_{p}\}$ are linearly independent. Let $~V_{u}^{j}$ be the subspace spanned by the 
 $\left\lbrace \nabla_{\omega,x}f_{k}(u,t),~k\neq j\right\rbrace .$\\
  we consider vectors
   $$N(f,j,u)=\nabla_{\omega,x}f_{j}(u,t)-q_{j}(u,t)~~1\leq j\leq p,$$
   where $q_{j}(u,t)$ is the projection of $\nabla_{\omega,x}f_{j}(u,t)$ to the subspace $V_{u}^{j}.$
 So there exist  $\alpha_{ji}$ such that 
  $$
  N(f,j,x)=\nabla_{\omega,x}f_{j}-\sum^{p}_{i=1,i\neq j}\alpha_{ji}\nabla_{\omega,x}f_{i}(u).
  $$
 
   When 
 $$
 u=(x,~\lambda)\in H_{\omega,r}^{\Sigma}(f,~\overline{w})\bigcap\{\|u\|<\alpha\}\backslash \Sigma\times \mathbf{R}^{l},
 $$
  by the condition $(K_{\Sigma}^{r,\delta})_{\Gamma},$ 
 if $\lambda_{k}(u,~t) \neq 0,$ 
 $$
 \frac{ \|\lambda_{k} \left( \nabla_{\omega,x}F_{k}(u,t)-\nabla_{\omega,x}f_{k}(u)\right) \|_{\omega}}{\Vert\sum^{p}_{j=1}\lambda_{j}\nabla_{\omega,x}f_{j}(u)\Vert_{\omega}}=
 \frac{ \| \nabla_{\omega,x}F_{k}(u,t)-\nabla_{\omega,x}f_{k}(u)\|_{\omega}}{\Vert\nabla_{\omega,x}f_{k}+\sum^{p}_{j=1,j\neq k}\left( \frac{\lambda_{j}}{\lambda_{k}}\right) \nabla_{\omega,x}f_{j}(u)\Vert_{\omega}}\leq \frac{o(d_{\omega}(x,~\Sigma)^{r})}{ C\cdot d_{\omega}(x,~\Sigma)^{r-\delta}}.
 $$ 
Hence, for a enough small positive number $\varepsilon,$
 \begin{eqnarray}
  \| \lambda_{k}\left( \nabla_{\omega,x}F_{k}(u,t)-\nabla_{\omega,x}f_{k}(u)\right) \|_{\omega}\leq \varepsilon d_{\omega}(x,~\Sigma)^{\delta}{\Vert\sum^{p}_{j=1}\lambda_{j}\nabla_{\omega,x}f_{j}(u)\Vert_{\omega}.}
    \end{eqnarray}

    Again since 
 \begin{eqnarray*}
  \Vert\sum^{p}_{j=1}\lambda_{j}\nabla_{\omega,x}F_{j}(u)\Vert_{\omega}&\geq & \Vert\sum^{p}_{j=1}\lambda_{j}\nabla_{\omega,x}f_{j}(u)\Vert_{\omega}-\Vert\sum^{p}_{j=1}\lambda_{j}\left(\nabla_{\omega,x}F_{j}-\nabla_{\omega,x}f_{j}\right) \Vert_{\omega}\\
&\geq & \Vert\sum^{p}_{j=1}\lambda_{j}\nabla_{\omega,x}f_{j}(u)\Vert_{\omega}-\sum^{p}_{j=1}\Vert\lambda_{j}\left(\nabla_{\omega,x}F_{j}-\nabla_{\omega,x}f_{j}\right) \Vert_{\omega}
 \end{eqnarray*}
 
So, by (3.9)
\begin{eqnarray*}
 \Vert\sum^{p}_{j=1}\lambda_{j}\nabla_{\omega,x}F_{j}(u)\Vert_{\omega}&\geq &  \Vert\sum^{p}_{j=1}\lambda_{j}\nabla_{\omega,x}f_{j}(u)\Vert_{\omega}-\sum^{p}_{j=1}\varepsilon d_{\omega}(x,~\Sigma)^{\delta}\Vert\sum^{p}_{j=1}\lambda_{j}\nabla_{\omega,x}f_{j}(u)\Vert_{\omega}\\
 &=&\left( 1-p\varepsilon d_{\omega}(x,~\Sigma)^{\delta}\right) \Vert\sum^{p}_{j=1}\lambda_{j}\nabla_{\omega,x}f_{j}(u)\Vert_{\omega}.
  \end{eqnarray*}
  Therefore
  \begin{eqnarray*}
 &\parallel & \nabla_{\omega,x}F_{k}-\sum^{p}_{j=1,j\neq k}\lambda_{j}\nabla_{\omega,x}F_{j}(u)\parallel_{\omega}\geq
 \left( 1-p\varepsilon d_{\omega}(x,~\Sigma)^{\delta}\right) \Vert\nabla_{\omega,x}f_{k}-\sum^{p}_{j=1,j\neq k}\lambda_{j}\nabla_{\omega,x}f_{j}(u)\Vert_{\omega}\\
&\geq &\left[ 1-\varepsilon\cdot p\cdot d_{\omega}(x,~\Sigma)^{\delta}\right] N(f,k,x)\geq\\
&\geq &  \left[ 1-\varepsilon\cdot p\cdot d_{\omega}(x,~\Sigma)^{\delta}\right]\cdot C\cdot d(x,~\Sigma)^{r-\delta}\geq\\
&\geq & C^{\prime}\cdot d_{\omega}(x,~\Sigma)^{r-\delta}
   \end{eqnarray*}
 is true for $ d_{\omega}(x,~\Sigma)<1.$
   I.e.
 $$
  d_{\omega,x}\nabla F=d_{\omega,x}(\nabla_{\omega,x}F_{1},\ldots,~\nabla_{\omega,x}F_{p})=
 min_{1\leq i\leq p}\lbrace\Vert N(F,i,x)\Vert_{\omega}\rbrace\geq C^{\prime} d_{\omega}(x,~\Sigma)^{r-\delta}
 $$ holds in a enough small neighborhood of 0.\\
 
 \ \ \textbf{Lemma 3.3.} 
{\sl Let $\Gamma$ be a compact Lie group acting orthogonally on space on space $(\mathbf{R}^{n}\times \mathbf{R}^{l},~\langle\cdot\rangle)$ and
 space $(\mathbf{R}^{p},~\langle\cdot\rangle)$ and $\Sigma$ be a germ of closed set of $\mathbf{R}^{n}$ such that $0\in \Sigma$ and  $\gamma\cdot\Sigma\subset \Sigma,~\forall \gamma \in \Gamma$. Suppose $f:(\mathbf{R}^{n}\times \mathbf{R}^{l},0) \longrightarrow (\mathbf{R}^{p},0)$ is $\Gamma-$equiviant.  
$H_{r}^{\Sigma}(f,~\overline{w})$ is a horn neighbourhood of $f$, then there exists a $C^{r}-\Gamma-$invariant function $\chi:(\mathbf{R}^{n}\times \mathbf{R}^{l},0) \longrightarrow (\mathbf{R}^{p},0),$ such that ,with $u\in(\mathbf{R}^{n}\times \mathbf{R}^{l}\setminus \lbrace 0\rbrace,$
$$~~0\leq \chi(u)\leq 1,~~~
 \chi(u)= \left\{
 \begin {array}{cc}
 1,& {\rm if}~ u\in H_{\omega,r}^{\Sigma}(f,~\frac{\overline{w}}{2}),
 \\[6pt]
 0,& {\rm if}~ u\in\mathbf{R}^{n}\times \mathbf{R}^{l}\setminus H_{\omega,r}^{\Sigma}(f,~\overline{w}). \\[6pt]
 \end {array}
        \right.
 $$
and $\chi(\gamma u)=\chi(u),~~\gamma\in\Gamma.$}\\

{\sl Proof.}~ Since $\Gamma\subset O(n),$ then, for $\gamma\in \Gamma,~\gamma\cdot H_{\omega,r}^{\Sigma}(f,~\overline{w})=H_{\omega,r}^{\Sigma}(f,~\overline{w})$ by Lemma 3.2(3). Let a function $\alpha:~\textbf{R}\longrightarrow \textbf{R}$ be defined by
$$
 \alpha(t)= \left\{
 \begin {array}{cc}
 ~e^{\frac{1}{t^{2}}},& {\rm if}~t~>0,
 \\[6pt]
~ 0,& {\rm if}~t\leq 0.\\[6pt]
 \end {array}
        \right.
 $$
$\alpha$ is a $C^{\infty}$ function. Again let a function $\beta:~\mathbf{R}^{n}\times \mathbf{R}^{l}\longrightarrow \textbf{R}$ be defined by

$$
\beta(x)=\frac{\alpha(1-\Vert x \Vert)}{\alpha(1-\Vert x \Vert)+\alpha(\Vert x \Vert-\frac{1}{2})}
$$
Owing to $\langle \gamma x,~\gamma x\rangle=\langle x,~x\rangle$ for~$~\gamma\in \Gamma,$ then $\beta(\gamma x)=\beta(x), ~\beta(x)$ is $C^{\infty}-\Gamma-$invariant function and $~0\leq \beta(x)\leq 1.$\\
\indent Moreover when $\Vert x \Vert\leq \frac{1}{2},~\beta(x)=1;$when $\Vert x \Vert>1,~\beta(x)=0.$\\
\indent Now the function $\chi:~(\mathbf{R}^{n}\times \mathbf{R}^{l},0) \longrightarrow \mathbf{R}$ is defined by $\chi(u)=\beta\left(\frac{f(u)}{\overline{w}\parallel u\parallel^{r}} \right) $
Then, for $u\in\mathbf{R}^{n}\times \mathbf{R}^{l}\setminus \lbrace 0\rbrace,$
$$~~0\leq \chi(u)\leq 1,~~~
 \chi(u)= \left\{
 \begin {array}{cc}
 ~~1,& {\rm if}~u\in H_{\omega,r}^{\Sigma}(f,~\frac{\overline{w}}{2}),
 \\[6pt]
~~ 0,& {\rm if}~u\in\mathbf{R}^{n}\times \mathbf{R}^{l}\setminus H_{\omega,r}^{\Sigma}(f,~\overline{w}). \\[6pt]
 \end {array}
        \right.
 $$
and $\chi(\gamma u)=\chi(u)=\beta\left(\frac{f(\gamma u)}{\overline{w}\parallel \gamma u\parallel^{r}} \right)=\chi(u)=\beta\left(\frac{\gamma f(u)}{\overline{w}\parallel u\parallel^{r}} \right)=\chi(u),~~\gamma\in\Gamma.$\\

  {\sl Proof of Theorem~1.1 }  Let $t_{0}$ be an arbitrary element of $I=[0,~1].$ Define
  $$F(u,~t)=f(u)+tp(u),~~u=(x,~\lambda))\in \mathbf{R}^{n}\times \mathbf{R}^{l},t\in I$$
where $p(u)$ satisfies
 \begin {eqnarray} 
\mid p_{i}\mid =o(d_{\omega}(x, \Sigma)^{r+\mid\omega\mid}),~~~~~~|\frac{\partial p_{i}}{\partial x_{j}}|=o(d_{\omega}(x, \Sigma)^{r}),~~~i=1,\ldots,p;~~j=1,\ldots,n.
\end {eqnarray} 
we define  in addition to the bilinear form on definition 1.1, 
$$
\langle\frac{\partial}{\partial u_{i}},~\frac{\partial}{\partial t}\rangle=0,~~1\leq i\leq n+l,~~\langle\frac{\partial}{\partial t},~\frac{\partial}{\partial t}\rangle=0.
$$

 \begin{eqnarray*}
  \nabla_{\omega}F_{i}& =&\sum^{n}_{j=1}\rho^{\omega_{j}}(\frac{\partial f_{i}}{\partial x_{j}}+t\frac{\partial p_{i}}{\partial x_{j}})\rho^{\omega_{j}}\frac{\partial }{\partial x_{j}}+\sum^{n+l}_{j=n+1}\rho^{\omega_{j}}(\frac{\partial f_{i}}{\partial \lambda_{j}}+t\frac{\partial p_{i}}{\partial \lambda_{j}})\rho^{\omega_{j}}\frac{\partial }{\partial\lambda_{j}}+p(u)\frac{\partial}{\partial t}\\
  &= & \nabla_{\omega,x}F_{i}+\sum^{n+l}_{j=n+1}\rho^{\omega_{i}}(\frac{\partial f_{i}}{\partial \lambda_{j}}+t\frac{\partial p_{i}}{\partial \lambda_{j}})\rho^{\omega_{j}}\frac{\partial }{\partial \lambda_{j}}+p(u)\frac{\partial}{\partial t}.
  \end{eqnarray*}
   we consider vectors
   $$N(F,j,u)=\nabla_{\omega,x}F_{j}(u,t)-Q_{j}(u,t)~~1\leq j\leq p,$$
   where $Q_{j}(u,t)$ is the projection of $\nabla_{\omega,x}F_{j}(u,t)$ to the subspace $W_{u}^{j},$ where $W_{u}^{j}$  is the subspace spanned by the 
 $\lbrace\nabla_{\omega,x}F_{k}(u,t);~k\neq j\rbrace.$\\
  \indent Since $f$ satisfies the  relative Kuo condition $(K_{\Sigma}^{r,\delta})_{\Gamma} $ and using Lemma 3.2,  there exist a positive number $ C^{\prime}$ such that,
  for
   $ (u,~t)\in W:= H_{r}^{\Sigma}(f,~\overline{w})\bigcap\{\|u\|<\alpha\})\times I, $
 $$
  d_{\omega,x}\nabla F= min_{1\leq i\leq p}\lbrace\Vert N(F,i,x)\Vert_{\omega}\rbrace\geq C^{\prime} d_{\omega}(x,~\Sigma)^{r-\delta}
 $$ holds. \\
 
\indent Firstly, we define a version the Kuo-vector field
  \[
 X_{1}(u,~t)= \left\{
 \begin {array}
 {r@{\quad,\quad}l}
 \frac{\partial}{\partial t}+\Sigma_{j=1}^{p}\frac{p_{i}(u)}{\Vert N(F,j,u)\Vert_{\omega}^{2}}N(F,j,u)&(u,~t)\in W\setminus\Sigma\times\mathbf{R}^{l}\times I\\
 
\frac{\partial}{\partial t}&(u,~t)\in W\cap \Sigma\times\mathbf{R}^{l}\times I \\[6pt]
 \end {array}
        \right.
\]
 \indent Now we give another form of $ X_{1}(u,~t)$.\\
   \indent In fact, since
 \begin{eqnarray*}
  \langle N(F,j,u),~\nabla_{\omega,x}F_{k}(u,t)\rangle &=&\langle\nabla_{\omega,x}F_{j}(u,t)-Q_{j}(u,t),\nabla_{\omega,x}F_{k}(u,t)\rangle\\
  &=&0~~1\leq j\leq p, ~k\neq j,
  \end{eqnarray*}
 where $Q_{j}(u,t)=\sum^{p}_{i=1,j\neq k}\alpha_{ij}\nabla_{\omega,x}F_{i}(u).$
 \indent Let $$
 b_{ki}=\langle \nabla_{\omega,x}F_{k}(u,t),~\nabla_{\omega,x}F_{i}(u,t)\rangle
 $$
 
  $$
   \left\{
 \begin {array}{ccc}

\sum^{p}_{k=1,j\neq k}b_{k1}\alpha_{kj}&+& b_{j1}\cdot(-1)=0\\ [6pt]

 \cdots &\cdots &\cdots\\  [6pt]
 \sum^{p}_{k=1,j\neq k}\hat{b_{kj} }\alpha_{kj}&+&\hat{b_{jj}}\cdot(-1)=0\\ [6pt]
  \cdots &\cdots &\cdots\\  [6pt]
\sum^{p}_{k=1,j\neq k}b_{kp}\alpha_{kj}&+&b_{jp}\cdot(-1)=0\\  [6pt]
\sum^{p}_{k=1,j\neq k}\nabla_{\omega,x}F_{k}(u,t)\alpha_{kj}&+ & Q_{j}(u,t)\cdot(-1)=0
\end{array} ,
 \right.$$
  where the hat means omission.\\
        
\indent This system of equations has non-zero solution $\alpha_{kj},\cdots,-1, j\neq k.$\\

 \indent So 
\begin{eqnarray*}
\left\vert 
\begin{array}{cccccc}
  b_{11}   &\cdots   &\hat{b_{1j}}  &\cdots & b_{1p}&  \nabla_{\omega,x}F_{1} \\
 \vdots & \vdots   &\vdots & \vdots & \vdots~& \vdots \\
   b_{j-1~ 1}   &\cdots   &\hat{b_{j-1~j}} &\cdots & b_{j-1~ p} &  \nabla_{\omega,x}F_{j-1} \\
    b_{j+1 ~1 }  &\cdots  &\hat{b_{j+1~j}} &\cdots & b_{j+1~p}&  \nabla_{\omega,x}F_{j+1} \\
    \vdots & \vdots   &\vdots & \vdots & \vdots~& \vdots\\
  b_{j1}  &\cdots   &\hat{b_{jj}}&\cdots  & b_{jp}&  Q_{j} \\
\end{array}
\right\vert=0.
\end{eqnarray*}

\indent  Now let $\left( D_{\omega,x}F\right)=\left(\nabla_{\omega,x}F_{1},\cdots,\nabla_{\omega,x}F_{p} \right)^{T} ,$
  $$
 \left( D_{\omega,x}F\right) \left( D_{\omega,x}F\right)^{T}=$$
 $$=
\left[ 
\begin{array}{cccc}
  \langle \nabla_{\omega,x}F_{1},~\nabla_{\omega,x}F_{1}\rangle   &\cdots   & \cdots &  \langle \nabla_{\omega,x}F_{1},~\nabla_{\omega,x}F_{p}\rangle \\
  \vdots   & \vdots  & \vdots~&  \vdots~\\
  
  \vdots   & \vdots & \vdots~&  \vdots~\\
  \langle \nabla_{\omega,x}F_{p},~\nabla_{\omega,x}F_{1}\rangle   &\cdots   & \cdots & \langle \nabla_{\omega,x}F_{p},~\nabla_{\omega,x}F_{p}\rangle .\\
\end{array}
\right] 
 $$
 $$=\left[ 
\begin{array}{cccc}
   b_{11}  &\cdots   & \cdots &   b_{1p} \\
  \vdots   & \vdots  & \vdots~&  \vdots~\\
  
  \vdots   & \vdots & \vdots~&  \vdots~\\
   b_{p1}   &\cdots   & \cdots &  b_{pp}.\\
\end{array}
\right] 
 $$
 If $M_{ij}$ denote the matrix obtained from $\left( D_{\omega,x}F\right) \left( D_{\omega,x}F\right)^{T}$
 deleting the $i$th row and the $j$th column. Let $A_{ij}=(-1)^{i+j}det(M_{ij}).$
then $Q_{j}$ equals to
 
  $$\frac{
 \left\vert
\begin{array}{cccccc}
  b_{11}   &\cdots   & b_{1j-1}  &  b_{1j+1} &\cdots &  \nabla_{\omega,x}F_{1} \\
 \vdots & \vdots   & \vdots & \vdots & \vdots~&\vdots\\
   b_{j-1~1}   &\cdots   & b_{j-1~j-1}  &  b_{j-1~j+1} &\cdots &  \nabla_{\omega,x}F_{j-1} \\
    b_{j+1~1}   &\cdots   & b_{j+1~j-1}  &  b_{j+1~j+1} &\cdots &  \nabla_{\omega,x}F_{j+1} \\
     \vdots & \vdots   & \vdots & \vdots & \vdots~&\vdots\\
   b_{j1}   &\cdots   & b_{j~j-1} &  b_{j~j+1} &\cdots &  0 \\
\end{array}
\right\vert}{-A_{jj}},
$$
i.e. $$
Q_{j}=-\sum_{i=1, i\neq j}^{p}\frac{A_{ji}}{A_{jj}}\nabla_{\omega,x}F_{i}.
$$
  Now $ N(F,j,u)=\nabla_{\omega,x}F_{j}(u,t)- Q_{j} $ is equal to
  $$
  \frac{
 \left\vert
\begin{array}{cccccc}
  b_{11}   &\cdots   & b_{1j-1}  &  b_{1j+1} &\cdots &  \nabla_{\omega,x}F_{1} \\
 \vdots & \vdots   & \vdots & \vdots & \vdots~&\vdots\\
   b_{j-1~1}   &\cdots   & b_{j-1~j-1}  &  b_{j-1~j+1} &\cdots &  \nabla_{\omega,x}F_{j-1} \\
    b_{j+1~1}   &\cdots   & b_{j+1~j-1}  &  b_{j+1~j+1} &\cdots &  \nabla_{\omega,x}F_{j+1} \\
     \vdots & \vdots   & \vdots & \vdots & \vdots~&\vdots\\
   b_{j1}   &\cdots   & b_{j~j-1} &  b_{j~j+1} &\cdots & \nabla_{\omega,x}F_{j}  \\
\end{array}
\right\vert}{-A_{jj}}
$$  
$~~~~=\sum_{i=1}^{p}\frac{A_{ji}}{A_{jj}}\nabla_{\omega,x}F_{i}.$ \\

 \indent  Owning  to 
  $$
  d_{\omega,x}\nabla F\geq C^{\prime} d(x,~\Sigma)^{r-\delta},~for ~(u,~t)=(x,~\lambda,~t)\in W\backslash \Sigma\times R^{l}\times I ,
 $$
  the vectors $\{\nabla_{\omega,x}F_{1},\ldots, \nabla_{\omega,x}F_{p}\}$ are linearly independent. So $ \left( D_{\omega,x}F\right) \left( D_{\omega,x}F\right)^{T}$
 has rank $p$ and
  $$
  \left( \left( D_{\omega,x}F\right) \left( D_{\omega,x}F\right)^{T}\right) ^{-1}=\frac{1}{d}
\left[ 
\begin{array}{cccc}
 A_{11}   &   A_{21}  &\cdots &  A_{p1} \\
  A_{12}   &   A_{22}  &\cdots &  A_{p2} \\
  \vdots  & \vdots    & ~\vdots & \vdots\\
   A_{1p}   &   A_{2p}  &\cdots &  A_{pp} \\
\end{array}
\right],
 $$
 where $d= \mid\left( D_{\omega,x}F\right) \left( D_{\omega,x}F\right)^{T}\mid.$

\indent We again compute $\| N(F,j,u)\|_{\omega}^{2}.$ 
 \begin{eqnarray*}
\| N(F,j,u)\|_{\omega}^{2}&=&\nabla_{\omega,x}F_{j}\cdot N(F,j,u)=\left\langle \nabla_{\omega,x}F_{j},~\sum_{i=1}^{p}\frac{A_{ji}}{A_{jj}}\nabla_{\omega,x}F_{i}\right\rangle\\
&=&\sum_{i=1}^{p}\frac{A_{ji}}{A_{jj}}\left\langle \nabla_{\omega,x}F_{j},\nabla_{\omega,x}F_{i}\right\rangle _{\omega}=\frac{d}{A_{jj}}.
 \end{eqnarray*}
 
So  \begin{eqnarray*}
\left(\frac{ N(F,1,u)}{\| N(F,1,u)\|_{\omega}^{2}},\cdots, \frac{ N(F,p,u)}{\| N(F,p,u)\|_{\omega}^{2}}\right)& =&\left(\frac{ N(F,1,u)}{\frac{d}{A_{11}}},\cdots, \frac{ N(F,p,u)}{\frac{d}{A_{pp}}}\right)\\
& =&\left(\frac{A_{11}}{d} N(F,1,u),\cdots,\frac{A_{pp}}{d} N(F,p,u)\right)
\end{eqnarray*}
Since
 $$
 \left( D_{\omega,x}F\right)^{T} \left( \left( D_{\omega,x}F\right) \left( D_{\omega,x}F\right)^{T}\right) ^{-1}=\left(\nabla_{\omega,x}F_{1},\cdots,\nabla_{\omega,x}F_{p} \right)\frac{1}{d}
\left[ 
\begin{array}{cccc}
 A_{11}   &   A_{21} &\cdots &  A_{p1} \\
  A_{12}   &   A_{22}  &\cdots &  A_{p2} \\
  \vdots  & \vdots    & \vdots~& \vdots\\
   A_{1p}   &   A_{2p}  &\cdots &  A_{pp} \\
\end{array}
\right]=
 $$
$$
=\left( \frac{1}{d}\left( \sum_{i=1}^{p}A_{1i}\nabla_{\omega,x}F_{i}\right) ,\cdots,\frac{1}{d}\left(\sum_{i=1}^{p}A_{pi}\nabla_{\omega,x}F_{i}\right) \right) 
 =\left(\frac{A_{11}}{d} N(F,1,u),\cdots,\frac{A_{pp}}{d} N(F,p,u)\right),
 $$
  Then,
 \begin{eqnarray}
\left(\frac{ N(F,1,u)}{\| N(F,1,u)\|_{\omega}^{2}},\cdots, \frac{ N(F,p,u)}{\| N(F,p,u)\|_{\omega}^{2}}\right)= \left( D_{\omega,x}F\right)^{T} \left( \left( D_{\omega,x}F\right) \left( D_{\omega,x}F\right)^{T}\right) ^{-1}.
\end{eqnarray}
  Therefore, when $~u\in H_{r}^{\Sigma}(f,~\overline{w})\bigcap\{\|u\|<\alpha\}\setminus\Sigma \times\mathbf{R}^{l} ,~t\in I,
   $
 \begin{eqnarray*}
 X_{1}(u,~t)&=&
  \frac{\partial}{\partial t}+\Sigma_{i=1}^{p}\frac{p_{i}(u)}{\Vert N_{i}\Vert_{\omega}^{2}}N_{i}\\
  ~~~~~& =&\left(-\Sigma_{i=1}^{p}\frac{p_{i}(u)}{\Vert N_{i}\Vert_{\omega}^{2}}N_{i} ,~0,~1 \right)^{T}\\
 ~~~~~ &=& \left(- \left[ \left( D_{\omega,x}F\right)^{T} \left( \left( D_{\omega,x}F\right) \left( D_{\omega,x}F\right)^{T}\right) ^{-1}p(u)\right]^{T} ,~0,~1 \right)^{T}.
 \end{eqnarray*}
 \indent Using the formula of $ X_{1}(u,~t),$  we show that $ X_{1}(u,~t)$ is $\Gamma-$equivariant.\\
 \indent Owning to Lemma 3.1 and $\Sigma$ is $\Gamma-$invariant. for $\forall\gamma \in\Gamma, H_{\omega,r}^{\Sigma}(f,~\overline{w} )$ is $\Gamma-$ invariant and $H_{\omega,r}^{\Sigma}(f,~\overline{w})\bigcap\{\|u\|<\alpha\}\setminus\Sigma \times\mathbf{R}^{l}$ is $\Gamma-$ invariant.  $~F(\gamma u,~t)=\gamma F( u,~t)$ with $\gamma\gamma^{T}=E.$
  Differentiating $ ~F_{i}(\gamma u,~t)=\gamma F_{i}( u,~t)$ with respect to $x$ yields
  \begin{eqnarray}
   \left(\frac{\partial ~F_{i}}{\partial x_{1}}(\gamma u,~t),\cdots,~\frac{\partial ~F_{i} }{\partial x_{n}}(\gamma u,~t)  \right)\gamma 
  =\gamma\left(\frac{\partial ~F_{i}( u,~t) }{\partial x_{1}},\cdots,~\frac{\partial ~F_{i}( u,~t) }{\partial x_{n}}  \right) .
   \end{eqnarray} 
   
By Lemma 3.1 and (3.11),  
\begin{eqnarray*}
\nabla_ {\omega,x}F_{i}(\gamma u,~t) &=&\left( \sum^{n}_{j=1}\rho(\gamma u)^{\omega_{j}}\frac{\partial F_{i}(\gamma u,~t)}{\partial x_{j}}\rho(\gamma u)^{\omega_{j}}\frac{\partial }{\partial x_{j}}\right) \gamma\\
&=&\gamma\left( \sum^{n}_{j=1}\rho( u)^{\omega_{j}}\frac{\partial F_{i}( u,~t)}{\partial x_{j}}\rho(u)^{\omega_{j}}\frac{\partial }{\partial x_{j}}\right). \\
\end{eqnarray*}
It implies $D_{\omega,x}F(\gamma u,~t)\gamma=\gamma D_{\omega,x} F( u,~t).$
  
 \begin{eqnarray*}
  &~& \left( D_{\omega,x}F(\gamma u,~t)\right)^{T} \left( \left( D_{\omega,x}F(\gamma u,~t)\right) \left( D_{\omega,x}F(\gamma u,~t)\right)^{T}\right) ^{-1}\\
  & =& \left(\gamma D_{\omega,x}F(u,~t)\gamma^{-1}\right)^{T} \left( \left(\gamma  D_{\omega,x}F(u,~t)\gamma^{-1}\right) \left( \gamma^{-1}\right)^{T}\left( D_{\omega,x}F(u,~t)\right) ^{T}\gamma^{T}\right) ^{-1}\\
  & =& (\gamma^{T})^{T}\left( D_{\omega,x}F(u,~t)\right)^{T}\gamma^{T} \left(\gamma  D_{\omega, x}F(u,~t)\gamma^{T}\gamma \left(D_{\omega,x}F(u,~t)\right)^{T}\gamma^{T}\right) ^{-1}\\
 & =& \gamma\left( D_{\omega,x}F(u,~t)\right)^{T}\gamma^{T}(\gamma^{T})^{T}  \left( D_{\omega,x}F(u,~t) \left( D_{\omega,x}F(u,~t)\right)^{T}\right) ^{-1} \gamma^{T}\\
 & =& \gamma\left( D_{\omega,x}F(u,~t)\right)^{T} \left( D_{\omega,x}F(u,~t) \left( D_{\omega,x}F(u,~t)\right)^{T}\right) ^{-1}\gamma^{-1}
 \end{eqnarray*}
 
    \begin{eqnarray*}
  X_{1}(\gamma u,~t)&=&  \left(- \left[ \left( D_{\omega,x}F(\gamma u,~t)\right)^{T} \left( \left( D_{\omega,x}F(\gamma u,~t)\right) \left( D_{\omega,x}F(\gamma u,~t)\right)^{T}\right) ^{-1}p(\gamma u)\right]^{T} ,~0,~1 \right)^{T}\\
  &=&\left(-\left[  \gamma\left( D_{\omega,x}F(u,~t)\right)^{T} \left( D_{\omega,x}F(u,~t) \left( D_{\omega,x}F(u,~t)\right)^{T}\right) ^{-1}\gamma^{-1}\gamma p(u)\right]^{T} ,~0,~1 \right)^{T}\\
   &=&\gamma\left(-\left[ \left( D_{\omega,x}F\right)^{T} \left( \left( D_{\omega,x}F\right) \left( D_{\omega,x}F\right)^{T}\right) ^{-1}p(u)\right]^{T} ,~0,~1 \right)^{T}\\
   &=&\gamma X_{1}( u,~t).
 \end{eqnarray*}
 
  \begin{eqnarray*}
  X_{1}(\gamma u,~t)&=&  \left(- \left[ \left( D_{\omega,x}F(\gamma u,~t)\right)^{T} \left( \left( D_{\omega,x}F(\gamma u,~t)\right) \left( D_{\omega,x}F(\gamma u,~t)\right)^{T}\right) ^{-1}p(\gamma u)\right]^{T} ,~0,~1 \right)^{T}\\
  &=&\left(-\left[  \gamma\left( D_{\omega,x}F(u,~t)\right)^{T} \left( D_{\omega,x}F(u,~t) \left( D_{\omega,x}F(u,~t)\right)^{T}\right) ^{-1}\gamma^{-1}\gamma p(u)\right]^{T} ,~0,~1 \right)^{T}\\
   &=&\gamma\left(-\left[ \left( D_{\omega,x}F\right)^{T} \left( \left( D_{\omega,x}F\right) \left( D_{\omega,x}F\right)^{T}\right) ^{-1}p(u)\right]^{T} ,~0,~1 \right)^{T}\\
   &=&\gamma X_{1}( u,~t).
 \end{eqnarray*}
 
  Therefore $ X_{1}(u,~t)$ is $\Gamma-$equivariant in W.\\
  \indent Again for 
  $~u\in H_{\omega,r}^{\Sigma}(f,~\overline{w})\bigcap\{\|u\|<\alpha\}\setminus\Sigma\times\mathbf{R}^{l} ,~t\in T,$
  \begin{eqnarray*}
 p_{j}(u) &-&\Sigma_{i=1}^{p}p_{i}(u)\nabla_{\omega,x}F_{j}\cdot\frac{ N(F,i,u)}{\Vert  N(F,i,u)\Vert_{\omega}^{2}}=\\
 &=& p_{j}(u) -\Sigma_{i=1}^{p}p_{i}(u)\langle\nabla_{\omega,x}F_{j}~,~N_{i}\rangle\cdot\frac{1}{\Vert N_{i}\Vert_{\omega}^{2}}=0,~~~j=1,\cdots,p.
 \end{eqnarray*}

 So
  \begin{eqnarray}
  \nabla_{\omega}F_{j}\cdot X_{1}(u,t)=<\nabla_{\omega}F_{j},~X_{1}(u,t)>=0, ~~~j=1,\cdots,p.
  \end{eqnarray} \\
 
\indent Moreover 
 $$
 <\Lambda,~X_{1}(u,t)>_{\omega}=0,~~\Lambda\in \lbrace 0\rbrace\times \mathbf{R}^{l}\times \lbrace 0\rbrace
 $$ 
 \begin{eqnarray}
 <\frac{\partial}{\partial t},~X_{1}(u,t)>=1
  \end{eqnarray}
 
and $ X_{1}(u,t)$ is $C^{1}$ on $ H_{r}^{\Sigma}(f,~\overline{w})\bigcap\{\|u\|<\alpha\})\setminus\Sigma\times\mathbf{R}^{l} ,~t\in T.
   $\\
   
 \indent By Lemma 3.3, we extend the vector field $~X_{1}(u,t)$ to a  vector field $~X(u,t)$ defined on a neighbourhood of zero.\\
 \indent Let  $~X(u,t)$ be defined by
 
\[
 X(u,~t)= \left\{
 \begin {array}
 {r@{\quad,\quad}l}
 \chi(u) X_{1}(u,~t)+ (1-\chi(u))\frac{\partial}{\partial t}& ~(u,t)\in\{\|u\|<\alpha\}\setminus\Sigma\times \mathbf{R}^{l}\times I
 \\[6pt]
\frac{\partial}{\partial t}&(u,t)\in\{\|u\|<\alpha\})\cap \Sigma\times \mathbf{R}^{l} \times I
 \end {array}
\right.\]

 \[
 ~~~~~~~~~~= \left\{
 \begin {array}
 {r@{\quad,\quad}l}
 \chi(u)\Sigma_{i=1}^{p}p_{i}(u)\frac{ N(F,i,u)}{\Vert  N(F,i,u)\Vert_{\omega}^{2}}+\frac{\partial}{\partial t}& (u,t)\in\{\|u\|<\alpha\}\setminus\Sigma\times \mathbf{R}^{l}\times I
 \\[6pt]
\frac{\partial}{\partial t}&(u,t)\in\{\|u\|<\alpha\})\cap \Sigma\times \mathbf{R}^{l} \times I
 \end {array}
\right.\]
i.e. when $ (u,t)\in\{\|u\|<\alpha\})\setminus \Sigma\times \mathbf{R}^{l} \times I,$
 $$
 X(u,~t)=\left(-\chi(u)\cdot\left[ \left( D_{\omega,x}F\right)^{T} \left( \left( D_{\omega,x}F\right) \left( D_{\omega,x}F\right)^{T}\right) ^{-1}p(u)\right]^{T} ,~0,~1 \right)^{T};
 $$
  when $(u,t)\in\{\|u\|<\alpha\})\cap \Sigma\times \mathbf{R}^{l} \times I,$
$$
 X(u,~t)=(0,~0,~1).
 $$
 Since
  \begin{eqnarray*}
\parallel X_{1}(u,~t) -\frac{\partial}{\partial t}\parallel_{\omega}&=&
  \parallel\Sigma_{j=1}^{p}\frac{p_{i}(u)}{\Vert N_{i}\Vert_{\omega}^{2}}N_{i}\parallel_{\omega}\\
 ~~~~~~~~~~~~ &\leq &\Sigma_{j=1}^{p}\frac{\mid p_{i}(u)\mid}{\Vert N_{i}\Vert_{\omega}} \\
  ~~~~~~~~~~~~  &\leq & \Sigma_{j=1}^{p}\frac{o\left( d_{\omega}(x,~\Sigma)^{r+\mid\omega\mid}\right) }{C^{'}d_{\omega}(x,~\Sigma)^{r-\delta}} \\
  ~~~~~~~~~~~~  &\leq & \Sigma_{j=1}^{p}\frac{o\left( d_{\omega}(x,~\Sigma)^{r+\mid\omega\mid}\right) }{C^{'}d_{\omega}(x,~\Sigma)^{r+\mid\omega\mid}}d_{\omega}(x,~\Sigma)^{r+\mid\omega\mid}d_{\omega}(x,~\Sigma)^{-r+\delta} \\ 
   ~~~~~~~~~~~~  &\leq & \Sigma_{j=1}^{p}\frac{o\left( d_{\omega}(x,~\Sigma)^{r+\mid\omega\mid}\right) }{C^{'}d_{\omega}(x,~\Sigma)^{r+\mid\omega\mid}}d_{\omega}(x,~\Sigma)^{\mid\omega\mid+\delta} ,\\ 
 \end{eqnarray*}

 then, by assumption $\mid\omega\mid+\delta-1>0,$
 \begin{eqnarray*}
 \parallel X(u,~t) -\frac{\partial}{\partial t}\parallel_{\omega}&= &
   \parallel \chi(u)\left(  X_{1}(u,~t)-\frac{\partial}{\partial t}\right)\parallel_{\omega}\\
   &\leq & \Sigma_{j=1}^{p}\frac{o\left( d_{\omega}(x,~\Sigma)^{r+\mid\omega\mid}\right) }{C^{'}d_{\omega}(x,~\Sigma)^{r+\mid\omega\mid}}d_{\omega}(x,~\Sigma)^{\mid\omega\mid+\delta-1}d_{\omega}(x,~\Sigma)\\
   &\leq & \Sigma_{j=1}^{p}\frac{o\left( d_{\omega}(x,~\Sigma)^{r+\mid\omega\mid}\right) }{C^{'}d_{\omega}(x,~\Sigma)^{r+\mid\omega\mid}}d_{\omega}(x,~\Sigma)\leq M d_{\omega}(x,~\Sigma),\\
    \end{eqnarray*}
  
 where $~M $ is a positive number and  $~u\in\{\|u\|<\alpha\})\setminus\Sigma\times \mathbf{R}^{l}$ with an enough small $\alpha$  and $d(x,~\Sigma)\leq 1,$\\
i.e. 
  \begin{eqnarray}
  \parallel X(u,~t) -\frac{\partial}{\partial t}\parallel_{\omega}\leq M d_{\omega}(x,~\Sigma),~~~for~~u\in\{\|u\|<\alpha\})\setminus\Sigma\times \mathbf{R}^{l}.
  \end{eqnarray} 
  
\indent Now we show $ X(u,~t)$ is $\Gamma-$equivariant. \\
\indent In fact,
\[
  X(\gamma u,~t)
= \left\{
  \begin {array}
 {r@{\quad,\quad}l}
 \chi(\gamma u) X_{1}(\gamma u,~t)+ (1-\chi(\gamma u))\frac{\partial}{\partial t}&~u\in\{\|u\|<\alpha\}\setminus\Sigma\times \mathbf{R}^{l} \\
\frac{\partial}{\partial t}&~u\in\{\|u\|<\alpha\}\cap \Sigma\times \mathbf{R}^{l} 
 \end {array}
        \right.\]

 \[ 
 ~ = \left\{ \begin {array}
 {r@{\quad,\quad}l}
  \chi( u)\gamma X_{1}(u,~t)+ (1-\chi( u))\frac{\partial}{\partial t}&~u\in\{\|u\|<\alpha\}\setminus\Sigma\times \mathbf{R}^{l}
 \\[6pt]
\frac{\partial}{\partial t}&~u\in\{\|u\|<\alpha\}\cap \Sigma\times \mathbf{R}^{l} 
 \end {array}
  \right.\]
 
  \[ 
 ~ = \left\{ \begin {array}
 {r@{\quad,\quad}l}
   \left(-\chi(u)\cdot\gamma\cdot\left[ \left( D_{x}F\right)^{T} \left( \left( D_{x}F\right) \left( D_{x}F\right)^{T}\right) ^{-1}p(u)\right]^{T} ,~0,~1 \right)^{T} &~u\in\{\|u\|<\alpha\}\setminus\Sigma\times \mathbf{R}^{l} \\
(0,~0,~1)&~u\in\{\|u\|<\alpha\}\cap \Sigma\times \mathbf{R}^{l} 
 \end {array}
        \right.
  \]
  
   \[ 
 ~ = \left\{ \begin {array}
 {r@{\quad,\quad}l}
  \gamma\cdot \left(-\chi(u)\cdot\left[ \left( D_{x}F\right)^{T} \left( \left( D_{x}F\right) \left( D_{x}F\right)^{T}\right) ^{-1}p(u)\right]^{T} ,~0,~1 \right)^{T} &~u\in\{\|u\|<\alpha\}\setminus\Sigma\times \mathbf{R}^{l}
 \\
(0,~0,~1)&~u\in\{\|u\|<\alpha\}\cap \Sigma\times \mathbf{R}^{l} 
 \end {array}
        \right.
  \]
  $=\gamma\cdot  X( u,~t).$\\
  
\indent Therefore $ X(u,~t)$ is $\Gamma-$equivariant.\\
\indent Moreover   $ X(u,~t)$ satisfies 
\begin{eqnarray}
 <\Lambda,~X(u,t)>=0,~~\Lambda\in \lbrace 0\rbrace\times \mathbf{R}^{l}\times \lbrace 0\rbrace
\end{eqnarray}
 
 \begin{eqnarray}
  <\frac{\partial}{\partial t},~X(u,t)>=0 
  \end{eqnarray}
 $ X(u,t)$ is $C^{1}$ on $\{\|u\|<\alpha\})\setminus\Sigma\times\mathbf{R}^{l}=\{\|u\|<\alpha\}\setminus\Sigma\times\mathbf{R}^{l} ,~t\in T.
  $ and 
 satisfies 
  \begin{eqnarray}
 \nabla_{\omega}F_{j}\cdot X(u,t)=<\nabla_{\omega}F_{j},~X(u,t)>=0, ~~~j=1,\cdots,p.
   \end{eqnarray}\\
  It implies $ X(u,t)$ is perpendicular to $ \nabla_{\omega,x}F$ at every $(u,~t)\in\{\|u\|<\alpha\})\setminus\Sigma\times\mathbf{R}^{l}\times I $, hence  $ X(u,t)$ is tangent to the level surface $F=$constant.\\
\indent In addition, if, for a enough small $\alpha,$ 
 $$
 u\in F^{-1}(0)\cap\{\|u\|<\alpha\}\setminus\Sigma\times\mathbf{R}^{l},
 $$
 then  $$
 ~u\in H_{\omega,r}^{\Sigma}(f,~w)\bigcap\{\|u\|<\alpha\}\setminus\Sigma\times\mathbf{R}^{l}.
 $$ 
\indent In fact, by $F_{i}(u,~t)=f_{i}(u)+tp_{i}(u)=0,$
 \begin{eqnarray*}
\frac{\mid f_{i}(u)\mid}{d_{\omega}(x,~\Sigma)^{r}}&=&\frac{\mid t p_{i}(u)\mid}{d_{\omega}(x,~\Sigma)^{r}}\\
&=&\frac{t\cdot o\left( d_{\omega}(x,~\Sigma)^{r+\mid\omega\mid}\right) }{d_{\omega}(x,~\Sigma)^{r}}\\
&<&\frac{ o\left( d_{\omega}(x,~\Sigma)^{r}\right) }{d_{\omega}(x,~\Sigma)^{r}}
 \end{eqnarray*}
 So when $\alpha$ is enough small and  
 $$
 u\in F^{-1}(0)\cap\{\|u\|<\alpha\}\setminus\Sigma\times\mathbf{R}^{l},
 $$ 
there is a enough small $\bar{w}$ such that  
$$
\parallel f(u)\parallel\leq \bar{w}  d_{\omega}(x,~\Sigma)^{r}\leq\bar{w}\| u\|_{\omega}^{r}=\bar{w}\rho(u)^{r}.
$$.
 \indent Therefore when $\alpha$ is enough small and  
 $$
 u\in F^{-1}(0)\cap\{\|u\|<\alpha\}\setminus\Sigma\times\mathbf{R}^{l},
 $$ we have 
 $$
 ~u\in H_{r}^{\Sigma}(f,~w)\bigcap\{\|u\|<\alpha\}\setminus\Sigma\times\mathbf{R}^{l}.
 $$   
 
\indent In order to be able to define $\phi$ , 
for a sufficient small $\alpha,$ if $(u,~t)\in \{\|u\|<\alpha\}\times I$, since (3.15), by                                                                                   Proposition 2.6, the following system of differential equations:
\begin{eqnarray}
u^{'}=X(u,~t)
\end{eqnarray}
is integrable.\\
\indent Now for $(u,~t)\in W$
 define $\gamma_{(u,t)}$ to be the maximal solution of $(3.19)$ such that $\gamma_{(u,~t)}(t)=u.$ 
 Let
 $H_{0},~\widetilde{H}_{0}:\{\|u\|<\alpha\}\times T\longrightarrow \{\|u\|<\alpha\}$
 be given by
\[~~~~~~~~~~H_{0}(u,~t)=\gamma_{(u,~t_{0})}(t),~~~~~~~~~~~~~~~~~~~\widetilde{H}_{0}(v,~t)=\gamma_{(v,t)}(t_{0}),\]
where $T$ is a small neighbourhood of $t$ in I.
By Proposition 2.6, the mappings $H_{0},\widetilde{H}_{0}$ are continuous mappings and uniqueness solutions of $(3.19)$ and the property $\gamma_{(\xi+\eta,x)}(s+t)=\gamma_{(\eta,\gamma_{(\xi,x)}(s))}(t)$ of the flow, it is easy to check that for $(u,~t)\in \{\|u\|<\alpha\}\times T$ we have
$$\widetilde{H}_{0}(H_{0}(u,t),t)=u,~H_{0}(u,t_{0})=u,$$ and 
$$H_{0}(\widetilde{H}_{0}(v,~t),t)=v,,\widetilde{H}_{0}(u,t)=u.$$ 
 $F(\gamma_{(u,t_{0}}(t),~t)=F(u,~t_{0})$ for any $t\in T,$
namely we have $f(~H_{0}(u,t))+t\cdot p(~H_{0}(u,t))=F(x,~t_{0})$ for $(u,~t)\in \{\|u\|<\alpha\}\times T.$ In particular, by (3.17),(3.18) for all $t,~t^{'}\in T,$ The germ of $F(u,~t)=0$ and $F(u,~t^{'})=0$ are $\Sigma-$homeomorphic(i.e. by a homeomorphism in $\mathcal{R}_{\Sigma}^{fix}).$  By  $X(u,~t)$ is $\Gamma-$equivariant, the  $\Sigma-$homeomorphic between the germ of $F(u,~t)=0$ and $F(u,~t^{'})=0$ is $\Gamma-$equivariant.\\
\indent  Using compactness of [0, 1], we obtain that a homeomorphism $\psi$ which has the form 
$$\psi(u,~t)=\left(\bar{\psi}(u,~t),~t \right)\in (\mathbf{R}^{n}\times \mathbf{R}^{l})\times I .$$
We define $\varphi=\bar{\psi}(u,~1)$ which is $\Gamma-\Sigma-$ BD equivalence between $G=F(\cdot,~0)$ and $\widetilde{G}=F(\cdot,~1).$\\
\indent Now we need to upgrading the equivalence between $G$ and $\widetilde{G}$ to a contact equivalence.  
we construct a matrix valued map $\tau$ which completes a contact equivalence in the same way as Theorem 3.1 of [6].\\
\indent 
Let 
$$
P_{\phi}(u)=\widetilde{G}(\phi(u))-G(u)~~~for~ u\in  \{\|u\|<\alpha\}
$$
$ P_{\phi}(u)$  is $\Gamma-$equivariant.\\
\indent We shall need to the fact  that 
\begin{eqnarray}
 \parallel d_{\omega}(x,~\Sigma)^{-r}P_{\phi}(u)\parallel=o(1).
 \end{eqnarray}

\indent Let
$$
Q_{\psi}(u,~s)=F\left(\psi (u,~s)\right)-G(u) 
$$
and
$$
\sigma=d_{\omega} \left( \bar{\psi}(u,~s),~\Sigma\times \mathbf{R}^{l}\right) .
$$
Then

\begin{eqnarray*}
 \frac{\partial}{\partial s}[Q_{\psi}(u,~s)]&=&F^{'}\left(\psi (u,~s)\right)\cdot X\left(\psi (u,~s)\right)\\
  &=&\left(\nabla_{\omega, x}F,~\nabla_{\omega,\lambda}F,~p \right)\cdot X\left(\psi (u,~s)\right)\\ 
 &=&[1-\chi\left(\bar{\psi}(u,~s)\right) ]p\left( \bar{\psi}(u,~s)\right)
 \end{eqnarray*}
 by (3.11).\\
So, since $Q_{\psi}(u,~0)=0,$
\begin{eqnarray*}
\parallel d_{\omega}(x,~\Sigma)^{-r}Q_{\psi}(u,~t) \parallel &\leq &\int_{0}^{t}\|d_{\omega}(x,~\Sigma)^{-r} \frac{\partial}{\partial s}[Q_{\psi}(u,~s)]\|ds\\
~~~~~~~~~~~~~~~~~~~~~~~~~~~~~~~~~~~~~~~~~~~~~~~~~~~~ &\leq &\int_{0}^{t}\|d_{\omega}(x,~\Sigma)^{-r}P\left( \bar{\psi}(u,~s)\right)\|ds\\
 ~~~~~~~~~~~~~~~~~~~~~~~~~~~~~~~~~~~~~~~~~~~~~~~~~~~~&\leq &\int_{0}^{t}\mid\left( \sigma/ d_{\omega}(x,~\Sigma)\right)^{r}\mid\cdot\| \sigma^{-r}P\left( \bar{\psi}(u,~s)\right)\|ds.
\end{eqnarray*}
\indent Because $d_{\omega}\left(u,~ \Sigma\times \mathbf{R}^{l}\right)=d_{\omega}\left(x,~\Sigma\right), $ by Lemma 2.7, we have

\begin {eqnarray} 
 d_{\omega} (x,~\Sigma)e^{-C\mid t \mid}\leq d_{\omega}\left( \bar{\psi}(u,~s),~\Sigma\times \mathbf{R}^{l}\right)\leq d_{\omega} (x,~\Sigma)e^{C\mid t \mid}
 \end {eqnarray}
for $t\in (\alpha,~\beta). $

By (2.1) and (3.21),
 \begin{eqnarray}
  \parallel d_{\omega}(x,~\Sigma)^{-r}Q_{\psi}(u,~t) \parallel=o(1)~~~~
uniformly~ in~ s\in [0, 1].~
\end{eqnarray}

\indent Because $Q_{\psi}(u,~1)=P_{\phi}(u),$ this proves (3.20).

\indent Again owing to 
$$
F(\psi(u,~s)=G\left( \bar{\psi}(u,~s)\right)+sP\left( \bar{\psi}(u,~s)\right),
$$
using (2.1),(3.21) and (3.22), if $u\in H_{\omega,r}^{\Sigma}(G,\beta)\bigcap \{\parallel u\parallel<\alpha\},$ then

$$
\Vert \sigma^{-r} G\left( \bar{\psi}(u,~s)\right)\Vert = \Vert \sigma^{-r}\left[ F\left( \psi(u,~s)\right)-s P\left( \bar{\psi}(u,~s)\right) \right] \Vert\leq
$$
$$
  \leq   \left( d_{\omega}(x,~\Sigma)/\sigma \right)^{r}\cdot \left\lbrace d_{\omega}(x,~\Sigma)^{-r}\Vert G(u)\Vert+  d_{\omega}(x,~\Sigma)^{-r}\| Q_{\psi}(u,~s)\|\right\rbrace
+\|\sigma^{-r}P\left( \bar{\psi}(u,~s)\right)\| 
$$
$$
=\left( d_{\omega}(x,~\Sigma)/\sigma\right)^{r}\cdot d_{\omega}(x,~\Sigma)^{-r}\|G(u)\|+o(1)\leq  e^{cr}\beta+0(1)
$$
 uniformly for $s\in [0,~1].$ Therefore, we may choose $V\subset \{\|u\|<\alpha\} $
to be a neighbourhood of the origin small enough and  $\beta$ sufficiently small such that 
 $$
\psi(u,~s)\in \left(H_{\omega,r}^{\Sigma}(G,~\frac{\bar{w}}{2}) \right)\cap \{\|u\|<\alpha\}\times T~~~when~~
(u,~s)\in \left(H_{\omega,r}^{\Sigma}(G,~\beta) \right)\cap V \times T.
$$
Hence, by (3.18), $F$ is constant on the flow of $X(u,~t)$ which remain in $
\psi(u,~s)\in \left(H_{r}^{\Sigma}(G,~\frac{\bar{w}}{2}) \right)\cap \{\|u\|<\alpha\}\times T~,$ so,
when $u\in \left(H_{\omega,r}^{\Sigma}(G,~\beta) \right)\cap V ,$
\begin {eqnarray}
P_{\phi}(u)=F(\psi(u,~1))-F(\psi(u,~0))=0.
\end {eqnarray}

 \indent Finally we construct a $\Gamma-$equivariant matrix valued map $\tau:~~~(V, 0)\longrightarrow(\mathbf{R}^{n}\times \mathbf{R}^{l},0)$ which completes a contact equivalence
For $u\in V $ and $v\in \mathbf{R}^{p},$ let
$$
\theta:~V\longrightarrow\mathcal{L}(\mathbf{R}^{p},~~\mathbf{R}^{p})
$$
be defined by 
 \begin{eqnarray}
\theta(u)(v)=\left\{
 \begin {array}{cc}
  \frac{\langle G(u)~,~v\rangle}{\Vert G(u)\Vert^{2}} P_{\phi}(u),&~G(u)\neq 0,
 \\
0,&~G(u)= 0.
 \end {array}
        \right.
 \end{eqnarray}
 
For $\gamma\in \Gamma,$ 
\begin{eqnarray*}
\theta(\gamma u)(v)&=&\left\{
 \begin {array}{cc}
  \frac{\langle G(\gamma u),~v\rangle}{\Vert G(\gamma u)\Vert^{2}} P_{\phi}(\gamma u),&~G(\gamma u)\neq 0, \\
0,&~G(\gamma u)= 0.
 \end {array}
    \right.
\end{eqnarray*}
\begin{eqnarray*}
~~~~~~~~~~~~~~~~~~&=&\left\{
 \begin {array}{cc}
  \frac{\langle\gamma G(u),~\gamma ~\gamma^{t}~v\rangle}{\Vert G(u)\Vert^{2}}\gamma P_{\phi}(u),~~G(\gamma u)\neq 0
 \\[6pt]
0,~~~~~~G(\gamma u)= 0\\[6pt]
 \end {array}
  \right.
 \end{eqnarray*}
        
\begin{eqnarray*}~~~~~~~~~~~~~~~~~&=&\left\{
 \begin {array}{cc}
  \frac{\langle G(u),~\gamma^{t}~v\rangle}{\Vert G(u)\Vert^{2}}\gamma P_{\phi}(u),&~G(\gamma u)\neq 0,
 \\
0,&~G(\gamma u)= 0
 \end {array}
   \right.
\end{eqnarray*}

Meantime, 
\begin{eqnarray*}
\gamma\cdot \theta(\gamma u)\cdot\gamma^{-1}(v)&=&\gamma\cdot \theta(\gamma u)\cdot\gamma^{t}(v)
=\gamma\cdot \theta(\gamma u)\cdot(\gamma^{t}v)\\
~~~~~~~~~~~~~~~~~~~~~~&=&\left\{
 \begin {array}{cc}
  \frac{\langle G(u),~\gamma^{t}~v\rangle}{\Vert G(u)\Vert^{2}}\gamma P_{\phi}(u),~~G(\gamma u)\neq 0
 \\[6pt]
0,~~~~~~G(\gamma u)= 0\\[6pt]
 \end {array}
        \right.
 \end{eqnarray*}
So $$
\theta(\gamma u)=\gamma\cdot \theta(\gamma u)\cdot\gamma^{t}=\gamma\cdot \theta(\gamma u)\cdot\gamma^{-1},  
$$
i.e. $\theta( u)$ is $\Gamma-$matrix valued map.\\
\indent When $G(u)\neq 0$ and $P_{\phi}(u))\neq 0,~\theta( u)$ is just a rank one linear transformation designed so that
$$
\theta( u)\left(G(u) \right) =P_{\phi}(u)).
$$

By (3.23) and (3. 24), $\theta( u)=0$ for $u\in \left(H_{\omega,r}^{\Sigma}(G,~\beta) \right)\cap V, $ which is a neighbourhood of $[G\mid V\setminus\Sigma\times\mathbf{R}^{l} ]^{-1}(0),$ so $\theta$ is continuous on $V\setminus\Sigma\times\mathbf{R}^{l}.$
By (3.20) and (3. 24) ,
$$
\|\theta(u)\|=\|P_{\phi}(u)\|/\|G(u)\|\leq o(d_{\omega}(x,~\Sigma)^{r})/d_{\omega}(x,~\Sigma)^{r}=o(1),
$$  
for $u\in V\setminus \left(H_{r}^{\Sigma}(G,~\beta) \right)\cap V, $
so $\theta$ is continuous at every point of $(\Sigma\times\mathbf{R}^{l})\cap V.$ 
Because $\theta$ is continuous on $V$
and  $\theta(0)=0,$ let $V$ be enough small so that
$$
I+\theta(u)\in \mathcal{GL}(\mathbf{R}^{p}),~~~~u\in V.
$$
Let
$$
\tau(u)=\left( I+\theta(u)\right) ^{-1},~~~~~u\in V.
$$
Then $$
\tau:~~V\rightarrow \mathcal{GL}(\mathbf{R}^{p})
$$
is clearly continuous and 
$$
G(u)=\tau(u\left(\tilde{G}\left(\phi(u) \right)  \right) ),
$$
since 
$$
\left( I+\theta(u)\right) \left( G(u)\right)=G(u)+ P_{\phi}(u))=\tilde{G}\left(\phi(u) \right) .
$$\\

\textbf{ Remark 3.}If $\Gamma$ be a compact Lie group acting linearly on space on space $(\mathbf{R}^{n}\times \mathbf{R}^{l},~\langle\cdot\rangle)$ and
 space $(\mathbf{R}^{p},~\langle\cdot\rangle),$ then  we can define a new inner on $\mathbf{R}^{n}\times \mathbf{R}^{l}$ following([10]):
$$
\langle u,~v\rangle_{\Gamma}=\int_{\Gamma}\langle \gamma u,~\gamma v\rangle d\gamma,
$$
where $\int$ is Haar integral on $\Gamma,~u,~v\in\mathbf{R}^{n}\times \mathbf{R}^{l}.$  It is important that $\Gamma$ be acting orthogonally on new inner space $(\mathbf{R}^{n}\times \mathbf{R}^{l},~\langle\cdot\rangle_{\Gamma})$ and new inner space $(\mathbf{R}^{p},~\langle\cdot\rangle_{\Gamma}).$\\

\indent Now we show to what terms from the Taylor expansion at every point that belongs to a closed subset of $\mathbf{R}^{n}$ such that $0\in \Sigma$ may be omitted without changing the topological type determined by $G$ and the value of the bifurcation parameter $\lambda$.\\

  {\sl Proof of Theorem 1.2.} The proof will be similar to that given in Theorem 1.1. \\
    \indent Let
    $$
    F(u,~t)=f(u)+tp(u).
    $$
  \indent  In proof of Theorem 1.1, Let $\Gamma=\left\lbrace e \right\rbrace  $ and $\omega= \lbrace \omega_{1},\omega_{2},\cdots,\omega_{n+l}\rbrace= \lbrace 1,~1,\cdots,~1\rbrace$, then $d_{\omega}(x, \Sigma)$ be substituted by $ d(x, \Sigma).$  When $u=(x,~\lambda)\in H^{\Sigma}_{r}(f,~\bar{w})\cap \\\lbrace\parallel u\parallel<\alpha\rbrace $ and $f$ satisfies the condition $(K_{\Sigma}^{r,\delta}),$ we have $ d _{x}\nabla F\geq C^{'}d(x, \Sigma)^{r-\delta}.$ \\
  \indent Again we use a version the Kuo-vector field 
  $$
 X_{1}(u,~t)= \left\{
 \begin {array}{cc}
 \frac{\partial}{\partial t}+\Sigma_{j=1}^{p}\frac{p_{i}(u)}{\Vert N_{i}\Vert^{2}}N_{i},&~(u,~t)\in W\setminus\Sigma\times T,
 \\
\frac{\partial}{\partial t},&~(u,~t)\in W\cap \Sigma\times T .
 \end {array}
        \right.
 $$ 
 where $W=H^{\Sigma}_{r}(f,~\bar{w})\cap {\parallel u\parallel<\alpha}.$ \\
 \indent Moreover we have a vector field $X(u,~t).$\\
 \indent Finally, using Lemma 2.4, we  may obtain a  homeomorphism  between 
   $f$ and $f+p.$\\
   
  {\sl Proof of Corollary 1.3.} Let $h(u)=g(u)-f(u).$ Then $j^{r+2}h(u)=0$ for $u=(x,~\lambda)\in (\Sigma\times \mathbf{R}^{l})\cap U.$ By Lemma 2.8, 
$$
\parallel h(u)\parallel=o\left( \left( d(u,~\Sigma\times \mathbf{R}^{l})^{r+2}\right) \right) =o\left( \left( d(x,~\Sigma)^{r+2}\right) \right) .
$$ It implies 
$$
|h_{i}|=o\left( d(x, \Sigma)^{d+1}\right) ,~~~~~~|\frac{\partial h_{i}}{\partial x_{j}}|=o\left( d(x, \Sigma)^{d}\right) ,~~~i=1,\ldots,p;~~j=1,\ldots,n.
$$
From Theorem 1.1, we obtain that $f$ and $g$ is $\Sigma-C^{0}$-BD equivalent and $\Sigma-C^{0}$-contact equivalent.\\

  \section{{$C^{0-}$Finite  determination of the  bifurcation diagram}}
 \indent ~~~ In the section, we show that contact $\vert \nu\vert-$determination of the  bifurcation problem is a corollary of Theorem 3.1.\\
 \indent We begin by presenting notation and concepts needed from [6].\\
 \indent For $A\in \mathcal{L}(\mathbf{R}^{n}, ~\mathbf{R}^{p}), $ let 
 $$\kappa(A)=inf{\parallel \alpha^{t}A\parallel:~\alpha\in \mathbf{R}^{p}),~\parallel \alpha^{t}\parallel=1}.
  $$
  When $A\in \mathcal{L}(\mathbf{R}^{n}, ~\mathbf{R}^{p})$ and rank(A)=p, let 
  $$A^{+}=A^{t}(AA^{t})^{-1}.$$ 
  Obvious $AA^{+}=I$ and $\kappa(A)=\parallel A^{+}\parallel^{-1}$ by [5].
  \indent For $\rho>0$ and $\nu\in ~\mathbf{R}^{p},$ let
  $$\rho^{\nu}=diag(\rho^{\nu_{1}},\cdots,~\rho^{\nu_{p}});$$
  $$\rho^{1-\mid\nu\mid}=diag(\rho^{1-\mid\nu\mid},\cdots,~\rho^{1-\mid\nu\mid}).$$
  \indent When $ F:\left( \mathbf{R}^{n}\times\mathbf{R}^{l},~0 \right)\rightarrow \left(  \mathbf{R}^{p},~0\right)$ is a $C^{1}$ map and $\nu\in ~\mathbf{R}^{p},$ we say that $F$ is {\sl ND($\nu$)} if there exist $\varepsilon>0,~\delta>0$ and a neighbourhood $U$ of the origin in $\mathbf{R}^{n}\times\mathbf{R}^{l}$ for which, with $\rho=\|u\|,$ 
  \begin{eqnarray}
  \kappa\left( \rho^{1-\nu}\nabla_{x}F(u)\right)\geq\varepsilon~~~~~if~ u=(x,~\lambda)\in H_{\nu}(F,~\delta~)\bigcap U, 
  \end{eqnarray}
  where $$
  H_{\nu}(F,~\delta~)=\lbrace u\in \mathbf{R}^{n}\times\mathbf{R}^{l}:~\Vert\rho^{1-\nu}F(u)\Vert\leq\delta \rbrace,~~\nabla_{x}F(u)=\left(\nabla_{x}F_{1}(u),\cdots, \nabla_{x}F_{p}(u)\right). $$
  
   \textbf{Remak 4.} If F is ND($\nu$),then F is also ND($\mid\nu\mid$) by [6].\\
   
    \indent The following lemmas will be used to prove Theorem 1.4.\\
   
 \indent   \textbf{Lemma 4.1.} ([11], Lemma A.9 ){\sl  The operator norm $\parallel\cdot\parallel$ of a submatrix is bounded by one of the whole matrix. More precisely, if $A\in C^{m\times n}$ has the form
   
   $$
    A=\left (\begin{array}{c|c}
    A^{(1)}& A^{(2)}\\ \hline
    A^{(3)}& A^{(4)}
    \end{array}\right) 
   $$
   for matrices $A^{(l)}$, then  $\parallel A^{(l)}\parallel\leq \parallel A \parallel$ for $l=1,2,3,4.$
  In particular, any entry of $A$ satisfies $\mid A_{j,k}\mid\leq \parallel A \parallel.$ }\\
  
   \indent {\sl Proof.} We give the proof for  $A^{(l)}$. The other cases are analogous. Let  $A^{(l)}$ be of size $m_{1}\times n_{1}.$ Then for the vector $x^{(1)}\in C^{n_{1}},$ we have 
    \begin{eqnarray*}
    \parallel A^{(l)x^{(1)}}\parallel^{2}&\leq & \parallel A^{(l)x^{(1)}}\parallel^{2}+ \parallel A^{(3)x^{(1)}}\parallel^{2}\\
    &=&\parallel\left (\begin{array}{c}
    A^{(1)}\\
    A^{(3)}
    \end{array}\right)x^{(1)} 
    \parallel^{2}\\
     &=&\parallel A\left (\begin{array}{c}
    x^{(1)}\\
     O
    \end{array}\right) 
    \parallel^{2}.
     \end{eqnarray*}
 The set $T_{1}$ of vectors $\left (\begin{array}{c}
    x^{(1)}\\
     O
    \end{array}\right)\in C^{n}$ with $\parallel  x^{(1)} \parallel \leq 1$ is contained in the set $T=\lbrace x\in C^{n}:~\parallel  x \parallel \leq 1\rbrace.$ Therefore, the supremum over $ x^{(1)}\in T_{1}$ above is bounded by $sup_{ x\in T}\parallel A x \parallel^{2}=\parallel A\parallel. $ This concludes the proof.\\
    
    \indent \textbf{Lemma 4.2.} {\sl  Suppose $F:~\left( \mathbf{R}^{n}\times\mathbf{R}^{l},~0 \right)\rightarrow \left(  \mathbf{R}^{p},~0\right)$ is a $C^{1}$ map and $\nu=(\nu_{1},\cdots, ~\nu_{p})$ such that F is ND($\nu)$. Then 
   $$ d_{x}\nabla F\geq C d(x,~0)^{\mid \nu\mid-1}~when~ u=(x,~\lambda)\in H(F,~\mid\nu\mid,~\delta~)\bigcap U. $$}
   {\sl Proof.} By Remark 4, F is also ND($\mid\nu\mid$). Then 
 
   \begin{eqnarray*}
    \left( \rho^{1-\mid\nu\mid} \nabla_{x}F \right)^{+}&=&
 \left(\rho^{1-\mid\nu\mid} \nabla_{x}F \right)^{t}\left(\rho^{1-\mid\nu\mid}\nabla_{x}F\cdot\left( \rho^{1-\mid\nu\mid} \nabla_{x}F \right) ^{t}\right)^{-1} \\ 
  & = &\rho^{1-\mid\nu\mid}\left( \nabla_{x}F \right)^{t} \rho^{\mid\nu\mid-1} \left( \left( \nabla_{x}F\right) \cdot \left( \nabla_{x}F \right) ^{t}\right)^{-1} \rho^{\mid\nu\mid-1}\\ 
   & = &\left( \nabla_{x}F \right)^{t}  \left( \left( \nabla_{x}F\right) \cdot \left( \nabla_{x}F \right) ^{t}\right)^{-1} \rho^{\mid\nu\mid-1}\\
    & = &\left( \nabla_{x}F \right)^{+}\cdot\rho^{\mid\nu\mid-1}\\
     & = &\rho^{\mid\nu\mid-1}\left( \nabla_{x}F \right)^{+}
 \end{eqnarray*}
 
   By (3.11),
  $$
 \left( \nabla_{x}F \right)^{+}= \left(\frac{ N(F,1,u)}{\| N(F,1,u)\|},\cdots, \frac{ N(F,p,u)}{\| N(F,p,u)\|}\right)
  $$
 $$ 
  \left( \rho^{1-\mid\nu\mid} \nabla_{x}F \right)^{+} =\rho^{\mid\nu\mid-1}  \left(\frac{ N(F,1,u)}{\| N(F,1,u)\|^{2}},\cdots, \frac{ N(F,p,u)}{\| N(F,p,u)\|^{2}}\right)
  $$
  where $  \rho^{1-\mid\nu\mid}=diag\left(\rho^{1-\mid\nu\mid},\cdots,\rho^{1-\mid\nu\mid} \right),$ when 
  $\nu=1,~\rho^{1-\mid\nu\mid}=I.$\\
  
\indent For vectors $$
 \frac{\rho^{\mid\nu\mid-1} N(F,i,u)}{\| N(F,i,u)\|^{2}},~~i=1,\cdots,p
 $$, By Lemma 4.1, 
 \begin{eqnarray*}
  \sum_{j=1}^{n+l} \left( \frac{\rho^{\mid\nu\mid-1} N(F,i,u)_{j}}{\| N(F,i,u)\|^{2}}\right)^{2} & =&
  \frac{1}{\| N(F,i,u)\|^{4}} \sum_{j=1}^{n+l} \left( \rho^{\mid\nu\mid-1} N(F,i,u)_{j}\right)^{2} \\
  &\leq &\sum_{j=1}^{n+l}\parallel \left( \rho^{1-\mid\nu\mid} \nabla_{x}F \right)^{+}\parallel^{2},
  \end{eqnarray*}
  where $ N(F,i,u)_{j}$ is the j-component of  $ N(F,i,u),$
  i.e.
  $$
   \rho^{2(\mid\nu\mid-1)}\frac{1}{\| N(F,i,u)\|^{2}}\leq (n+l)\parallel \left( \rho^{1-\mid\nu\mid} \nabla_{x}F \right)^{+}\parallel^{2}.
  $$
  
  $$
  \kappa\left( \rho^{1-\mid\nu\mid}\nabla_{x}F(u)\right)=\parallel \left( \rho^{1-\mid\nu\mid} \nabla_{x}F \right)^{+}\parallel^{-1}\leq \left( \rho^{\mid\nu\mid-1}\frac{1}{\| N(F,i,u)\|}\right) ^{-1}\frac{1}{((n+l))^{\frac{1}{2}}},
  $$
 i.e.
$$
  \kappa\left( \rho^{1-\mid\nu\mid}\nabla_{x}F(u)\right)\leq \rho^{1-\mid\nu\mid}\| N(F,i,u)\|\frac{1}{\mid\nu\mid-1}.
 $$
  Since
   $$ d_{x}\nabla F= min \{\| N(F,i,u)\|:~~i=1,\cdots,~p\},$$ 
 
  and $$ 
  \kappa\left( \rho^{1-\mid\nu\mid}\nabla_{x}F(u)\right)\geq\varepsilon,
    $$
    then
   $$
   d_{x}\nabla F\geq\varepsilon((n+l))^{\frac{1}{2}}\rho^{\mid\nu\mid-1}\geq C d (x,0)^{\mid\nu\mid-1}
   $$
   when $u=(x,~\lambda)\in H_{\mid\nu\mid}(F,~\delta~)\bigcap U.$
   where $C=\varepsilon((n+l))^{\frac{1}{2}}.$\\
   
    {\sl Proof of Theorem 1. 4.} In Theorem 1.1, let $\Gamma={e},~\Sigma=\{0\}$ and $\omega=(1,1,\cdots,1).$
   Since F is 
   ND($\nu),$  F satisfies relative Kuo condition $\left( K_{\lbrace 0\rbrace}^{\mid\nu\mid,~1}\right)_{\lbrace e\rbrace} $ by Lemma 4.2. Again  Theorem 1.1, F is BD $\mid\nu\mid-$ determined.
\section*{Acknowledgements}
The authors thanks the referee for his/her careful reading and very
useful comments which improved the final version of this paper.


\begin{thebibliography}{0}
\bibitem{emot53 }K.Bekka, S. Koike, \emph{Characterisations of V-sufficiency and $C^{0}$-sufficiency of relative jets,} arXiv:1703.07069v4[math. AG] 17 Mar 2020.
\bibitem{emot53 }M. Buchner, J.Marsd$\acute{e}$n and S.Schecter, \emph{Applications of  the blowing -up construction and algebaic geometry to bifurcation problems,} J. Differential Equations, Vol.48,1983,404-433.

\bibitem{emot53}T. Fukui, L. Paunescu, \emph{Stratification theory from the weighted point of view}
Canad. J. Math. vol.53(1),(2001), 73-97.
\bibitem{c85} M.Golubitsky, D. Scheaffer, \emph{Singularities and groups in bifurcation theory, Vol.1}, Applied Mathematical Sciences 51. Spring-Verlag, 1985.
\bibitem{cmz96} M.Golubitsky,I. Stewart, D. Scheaffer, \emph{Singularities and groups in bifurcation theory, Vol.2}, Applied Mathematical Sciences 69. Spring-Verlag, 1988.
\bibitem{lz} P.B.Percell, P,N. Brown, \emph{Finite dermination of bifurcation problems},SIAM J. MATH. ANAL.  Vol.16, No.1,(1985),28-46.
\bibitem{hm89}L.Paunescu, \emph{A weighted version of the Kuiper-Kuo-Bochnack-Lojasiewics Theorem}, J. Algebraic Geometry, Vol.2, (1993), 66-79.
\bibitem{kr90}M.A.S.Ruas, M.J.Saia,\emph{ $C^{l}-$ determinacy of weighted homogeneous germs},
 Hokkaido Math.J., Vol.26,(1997), 89-99.
 \bibitem{kr90} T. Brocker, T. tom Dieck,\emph{ Representations of compact Lie groups}, GTM 98,
Spring-Verlag, New York, 1985.
 \bibitem{Z03}D. Bump, \emph{Lie groups, } GTM 225, Springer-Verlag, New York, 2004.
,1992 (Ch).
\bibitem{Z03} S.Foucart, H.Rauhut, \emph{A mathematical introduction to compressive sensing}, Springer New York, 2013.
\bibitem{Z03} L.Hengxing, Z.Dun-mu, \emph{$C^{l}-\mathbf{g}_{V}-$determinacy of weighted homogeneous function germs on weighted homogeneous analytic varieties}, Hokkaido Math.J. ,Vol.37, (2008), 309-329.
 \bibitem{Z03}
B. Osi$\acute{n}$ska-Ulrych, T. Rodak, G. Skalski, \emph{Topological triviality of deformations of regular mappings},  Bull. Sci.Math. Vol.161, (2020), 1-21.
 \bibitem{Z03}
 Z. Jiangcheng, S.Fuwei, S.Ruixia, L.Guofu, \emph{D-determination of bifurcation problems with respect to $C^{0}$ contact equivalence from the weighted point of view}, Acta. Math. Sinica (Chinese),Vol.42, No 2,(1999), 305-312.
 \bibitem{Z03}
 T.C.Kuo,\emph{ Characterizations of v-sufficiency of jets}, Topology, 11(1972), 115-131.
\end{thebibliography}
\end{document}